\documentclass{article}

\usepackage[utf8]{inputenc}
\usepackage[english]{babel}
\usepackage{xcolor}
\usepackage{hyperref}
\usepackage{cite}  
\usepackage{algorithm} 
\usepackage{algpseudocode} 
\usepackage{graphicx}
\usepackage{amsmath}
\usepackage{bm}
\usepackage{soul}
\usepackage{authblk}
\usepackage{url}
\usepackage{hyperref}

\usepackage[a4paper]{geometry}

\usepackage{algorithm}
\usepackage{algpseudocode}
\usepackage{amsfonts}
\usepackage{amsmath}
\usepackage{amssymb} 
\usepackage{amsthm} 
\usepackage{bm,bbm}
\usepackage{booktabs}
\usepackage{chapterbib} 
\usepackage{color}
\usepackage{comment}
\usepackage{enumitem}
\usepackage{epsfig}
\usepackage{graphicx,graphics,rotating}
\usepackage{multicol}
\usepackage{multirow}
\usepackage{overpic}
\usepackage{pgfplots}
\usepackage{pgfplotstable}
\pgfplotsset{compat=newest}
\usepackage{siunitx}
\usepackage{subcaption}
\usepackage{tikz}
\usepackage{times}
\usepackage{xcolor}

\theoremstyle{plain}
\newtheorem{theorem}{Theorem}[section]

\theoremstyle{definition}

\newtheorem{remark}[theorem]{Remark}

\theoremstyle{plain}

\newcommand{\uspace}{u_1}
\newcommand{\utime}{u_2}
\newcommand{\urank}{u_3}

\newcommand{\toll}{\epsilon}

\newcommand{\rate}{rate}
\newcommand{\cstrg}{strg}
\newcommand{\ctime}{time}

\def\trait #1 #2 #3 {\vrule width #1pt height #2pt depth #3pt}
\def\fin{\hfill
        \trait .3 5 0
        \trait 5 .3 0
        \kern-5pt
        \trait 5 5 -4.7
        \trait 0.3 5 0
\medskip}

\newif\ifANKI

\newcounter{numbs}
\newcounter{numbi}
\newcounter{numbii}

\ANKIfalse

\ifANKI

\else

\fi

\newcommand{\PARAGRAPH}[1]{\medskip\noindent\textit{\textbf{#1}}}

\newcommand{\REAL}{\mathbbm{R}}

\newcommand{\av}{\mathbf{a}}
\newcommand{\bv}{\mathbf{b}}
\newcommand{\ccv}{\mathbf{c}}

\newcommand{\rv}{\mathbf{r}}

\newcommand{\xv}{\mathbf{x}}

\newcommand{\as}{a}
\newcommand{\bs}{b}

\newcommand{\fs}{f}
\newcommand{\gs}{g}
\newcommand{\hs}{h}

\newcommand{\ms}{m}
\newcommand{\ns}{n}

\newcommand{\ps}{p}

\newcommand{\rs}{r}

\newcommand{\ts}{t}
\newcommand{\us}{u}

\newcommand{\xs}{x}
\newcommand{\ys}{y}
\newcommand{\zs}{z}

\newcommand{\As}{A}
\newcommand{\Bs}{B}
\newcommand{\Cs}{C}

\newcommand{\Is}{I}

\newcommand{\Ns}{N}

\newcommand{\Ps}{P}

\newcommand{\Ts}{T}
\newcommand{\Us}{U}

\newcommand{\Xs}{X}
\newcommand{\Ys}{Y}
\newcommand{\Zs}{Z}

\newcommand{\Nbnd}{\Ns^{\BND}}

\newcommand{\usb}{\overline{u}}

\newcommand{\ust}{\widetilde{u}}

\newcommand{\xst}{\widetilde{x}}
\newcommand{\yst}{\widetilde{y}}
\newcommand{\zst}{\widetilde{z}}

\newcommand{\ass}[1]{a_{#1}}
\newcommand{\bss}[1]{b_{#1}}

\newcommand{\iss}[1]{i_{#1}}
\newcommand{\jss}[1]{j_{#1}}
\newcommand{\kss}[1]{k_{#1}}

\newcommand{\nss}[1]{n_{#1}}

\newcommand{\rss}[1]{r_{#1}}

\newcommand{\uss}[1]{u_{#1}}

\newcommand{\xss}[1]{x_{#1}}
\newcommand{\yss}[1]{y_{#1}}
\newcommand{\zss}[1]{z_{#1}}

\newcommand{\Ass}[1]{A_{#1}}

\newcommand{\Pss}[1]{P_{#1}}

\newcommand{\Uss}[1]{U_{#1}}

\newcommand{\usz}[2]{\us_{#1}^{#2}}

\newcommand{\Usz}[2]{\Us_{#1}^{#2}}

\newcommand{\matD}{\mathsf{D}}

\newcommand{\matI}{\mathsf{I}}

\newcommand{\calI}{\mathcal{I}}

\newcommand{\calO}{\mathcal{O}}

\newcommand{\calR}{\mathcal{R}}

\newcommand{\tA}{\mathcal{A}}
\newcommand{\tB}{\mathcal{B}}
\newcommand{\tC}{\mathcal{C}}

\newcommand{\tF}{\mathcal{F}}

\newcommand{\tR}{\mathcal{R}}

\newcommand{\tU}{\mathcal{U}}

\newcommand{\tZ}{\mathcal{Z}}

\newcommand{\tATT}{\mathcal{A}^{\TT}}

\newcommand{\LTWO}  {L^2}

\newcommand{\CS}[1] {C^{#1}}

\renewcommand{\P} {\textsf{P}}            

\newcommand{\hh}{h}

\newcommand{\xvV}{\xv_{\V}}        
\newcommand{\xvC}{\xv_{\C}}        

\newcommand{\DIM} {d}              
\newcommand{\BND} {\text{bnd}}       

\newcommand{\NMB}{N}

\newcommand{\dt}{dt}

\newcommand{\psV}{\ps_{\V}}

\newcommand{\nlen}{\hspace{-0.2mm}}

\newcommand{\norm}   [2]{|\nlen|#1|\nlen|_{#2}}

\newcommand{\restrict}[2]{{#1}_{|{#2}}}

\newcommand{\EOD}{\end{document}}

\newcommand{\roundPrecision}{2}
\newcommand{\TT}{\scalebox{0.6}{{\it TT}}}

\newcommand{\NV}{\NMB_{\V}}
\newcommand{\NVx}{\NV^x}
\newcommand{\NVy}{\NV^y}
\newcommand{\NVz}{\NV^z}

\newcommand{\NC}{\NMB_{\C}}
\newcommand{\NCx}{\NC^x}
\newcommand{\NCy}{\NC^y}
\newcommand{\NCz}{\NC^z}

\newcommand{\NT}{\NMB_{\Delta\ts}}

\newcommand{\hx}{\hh_x}
\newcommand{\hy}{\hh_y}
\newcommand{\hz}{\hh_z}

\newcommand{\V} {\textit{v}} 
\newcommand{\C} {\textit{c}} 

\newcommand{\ic}{\iss{\C}}
\newcommand{\jc}{\jss{\C}}
\newcommand{\kc}{\kss{\C}}       

\newcommand{\iV}{\iss{\V}}
\newcommand{\jV}{\jss{\V}}
\newcommand{\kV}{\kss{\V}}       

\newcommand{\xsV}{\xss{\V}}
\newcommand{\ysV}{\yss{\V}}
\newcommand{\zsV}{\zss{\V}}
\newcommand{\usV}{\uss{\V}}

\newcommand{\xsC}{\xss{\C}}
\newcommand{\ysC}{\yss{\C}}
\newcommand{\zsC}{\zss{\C}}
\newcommand{\usC}{\uss{\C}}

\newcommand{\xsc}{\xs(\ic)}
\newcommand{\ysc}{\ys(\jc)}
\newcommand{\zsc}{\zs(\kc)}

\newcommand{\RESD}{\calR}
\newcommand{\RESDV}{\calR_{\V}}

\newcommand{\IntpV}{\calI^{\V}}

\newcommand{\tAss}[1]{\tA_{#1}}
\newcommand{\tRss}[1]{\tR_{#1}}
\newcommand{\tUss}[1]{\tU_{#1}}
\newcommand{\tZss}[1]{\tZ_{#1}}

\newcommand{\usVx}{\partial^{c}_{x}\usV}
\newcommand{\usVy}{\partial^{c}_{y}\usV}
\newcommand{\usVz}{\partial^{c}_{z}\usV}

\newcommand{\usVxx}{\partial^{c}_{xx}\usV}
\newcommand{\usVyy}{\partial^{c}_{yy}\usV}
\newcommand{\usVzz}{\partial^{c}_{zz}\usV}

\newcommand{\fsV} {\fs_{\V}}
\newcommand{\fsVTT} {\fsV^{\TT{}}}

\newcommand{\rsVTT}{\rsV^{\TT{}}}
\newcommand{\psVTT}{\psV^{\TT{}}}
\newcommand{\usVTT} {\usV^{\TT{}}}

\newcommand{\usVTTx}{\partial^{c}_{x}\usVTT}
\newcommand{\usVTTy}{\partial^{c}_{y}\usVTT}
\newcommand{\usVTTz}{\partial^{c}_{z}\usVTT}

\newcommand{\usVTTxx}{\partial^{c}_{xx}\usVTT}
\newcommand{\usVTTyy}{\partial^{c}_{yy}\usVTT}
\newcommand{\usVTTzz}{\partial^{c}_{zz}\usVTT}

\newcommand{\usVTTxy}{\partial^{c}_{xy}\usVTT}
\newcommand{\usVTTyz}{\partial^{c}_{yz}\usVTT}
\newcommand{\usVTTzx}{\partial^{c}_{zx}\usVTT}

\newcommand{\usVTTyx}{\partial^{c}_{yx}\usVTT}
\newcommand{\usVTTzy}{\partial^{c}_{zy}\usVTT}
\newcommand{\usVTTxz}{\partial^{c}_{xz}\usVTT}

\newcommand{\DeltaC}{\Delta^{c}}
\newcommand{\DeltaV}{\Delta^{\V}}

\newcommand{\DeltaCTT}{\Delta^{c,\TT{}}}
\newcommand{\DeltaVTT}{\Delta^{\V,\TT{}}}

\newcommand{\IntpVTT}{\calI^{\V,\TT{}}}

\newcommand{\hxss}[1]{\hh_{x,#1}}
\newcommand{\hyss}[1]{\hh_{y,#1}}
\newcommand{\hzss}[1]{\hh_{z,#1}}

\newcommand{\usVb}{\usb_{\V}}
\newcommand{\usVTTb}{\usb^{\TT}_{\V}}
\newcommand{\RNDG}{\scalebox{0.9}{{\sf rndg}}}
\newcommand{\TTSVD}{\scalebox{0.9}{{\sf TT-SVD}}}

\newcommand{\cTT}{\mathbbm{TT}}

\newcommand{\SFFONT}[1]{\scalebox{0.9}{{\sf #1}}}

\newcommand{\dirX}{\SFFONT{X}}
\newcommand{\dirY}{\SFFONT{Y}}
\newcommand{\dirZ}{\SFFONT{Z}}

\newcommand{\SET}[1]{\big(#1\big)}

\newcommand{\rsV}{\rs_{\V}}

\newcommand{\bsV}{\bs_{\V}}

\newcommand{\ztVTT}{\zeta^{\TT}_{\V}}
\newcommand{\bsVTT}{\bs^{\TT}_{\V}}

\newcommand{\MATVECPROD}{\scalebox{0.9}{{\sf Matrix\_Vector\_Product}}}
\newcommand{\PRECON}{\scalebox{0.9}{{\sf Preconditioner}}}
\newcommand{\MAXVOL}{\scalebox{0.9}{{\sf MaxVol}}}
\newcommand{\TRIDIAG}{\scalebox{0.9}{{\sf tridiag}}}

\begin{document}

  \title{The low-rank tensor-train finite difference method\\
    for three-dimensional parabolic equations}

  \author[1] {Gianmarco Manzini}
  \author[2] {Tommaso Sorgente}

  \affil[1]{T-5, Theoretical Division, Los Alamos National Laboratory, Los Alamos, NM, USA}
  
  \affil[2]{Istituto di Matematica Applicata e Tecnologie Informatiche ``Enrico~Magenes'', 
  Consiglio Nazionale delle Ricerche, Genova, Italy}
  
  \date{}

\maketitle

\begin{abstract}
This paper presents a numerical framework for the low-rank
approximation of the solution to three-dimensional parabolic problems.
The key contribution of this work is the tensorization process based
on a tensor-train reformulation of the second-order accurate finite
difference method.
We advance the solution in time by combining the finite difference method
with an explicit and implicit Euler method and with the Crank-Nicolson
method.
We solve the linear system arising at each time step from the implicit
and semi-implicit time-marching schemes through a matrix-free
preconditioned conjugate gradient (PCG) method, appositely designed to
exploit the separation of variables induced by the tensor-train format.
We assess the performance of our method through extensive numerical experimentation, demonstrating that the tensor-train design offers a robust and highly efficient alternative to the traditional approach. 
Indeed, the usage of this type of representation leads to massive time and memory savings while guaranteeing almost identical accuracy with respect to the traditional one.
These features make the method particularly suitable to tackle challenging high-dimensional problems.
\end{abstract}


\section{Introduction}
\label{sec:model}
We consider the three-dimensional, open, bounded, hyper-rectangular
domain $\Omega\subset\REAL^3$ with boundary $\Gamma=\partial\Omega$.
We are interested in the numerical approximation to the solution $\us$
of the parabolic equation with Dirichlet boundary condition:
\begin{equation}
  \begin{aligned}
    \frac{\partial\us}{\partial\ts} - \Delta\us &= \fs & \qquad\text{in }\Omega\times(0,\Ts], \\
    \us &= \gs                                         & \qquad\text{on }\Gamma\times(0,\Ts],
  \end{aligned}
  \label{eq:problem:strong}
\end{equation}
where $\fs$ and $\gs$ on the right-hand side of
\eqref{eq:problem:strong} are the forcing term and the Dirichlet
boundary function, respectively.
To complete the mathematical model, we assume that $\us$ at the
initial time $t=0$ satisfies the condition:
\begin{align}
  \us(\cdot,0)=\us_0 \qquad\text{in }\overline{\Omega},
  \label{eq:problem:strong2}
\end{align}
where $\overline{\Omega}$ is the closure of $\Omega$ in $\REAL^3$.
Under suitable regularity assumption on $f,g,u_0$, and $\Gamma$, we
can prove that the problem is well-posed~\cite{evans2010partial}.
Throughout this paper, we assume that the solution
to~\eqref{eq:problem:strong} exists, is unique, and possesses
sufficient regularity for all the approximations and estimates that
follow to be meaningful.

Numerous approximation methods have been developed to solve such model
problems, including
\emph{Finite Difference Methods (FDMs)}~\cite{Strikwerda:2007},
\emph{Finite Element Methods (FEMs)}~\cite{Braess:2007},
\emph{Finite Volume Methods (FVMs)}~\cite{Versteeg-Malalasekera:2007},
and
\emph{Spectral Methods}~\cite{Shen-Tang-Wang:2011,Fornberg:1996}.
Time integration typically employs explicit, implicit, or
semi-implicit (Crank-Nicolson) schemes~\cite{Leveque:2007}.
While implicit and semi-implicit methods offer better stability than
explicit schemes, allowing larger time steps, they require solving
linear systems at each time step, increasing the computational
complexity.

\subsection{A first step towards dual grid discretizations}

In this work, we present a second-order accurate spatial
discretization for
problem~\eqref{eq:problem:strong}-\eqref{eq:problem:strong2},
employing a \textit{dual grid approximation} of $\us$ and its
Laplacian, $\Delta\us$, which combines both cell-centered and
vertex-centered discretizations of these quantities.
The dual grid approach offers several decisive advantages over
traditional single grid methods that make it particularly suitable for
tensor-train reformulations.
First, dual grid discretizations eliminate spurious numerical
oscillations that commonly plague single grid methods in fluid
dynamics applications.
As demonstrated in the pioneering work by Harlow and
Welch~\cite{Harlow-Welch:1965}, positioning velocity components at
cell faces while maintaining pressure at cell centers creates natural
pressure-velocity coupling that prevents checkerboard pressure modes
without requiring artificial stabilization techniques.
Second, dual grid methods achieve exact discrete conservation of
fundamental physical quantities in fluid dynamics simulation, e.g.,
mass, momentum, and energy, at machine precision through their
inherent geometric
structure~\cite{Perot:2011,Morinishi-Lund-Vasilyev-Moin:1998}.
This property, impossible to guarantee with interpolation-based single
grid approaches, prevents unphysical solution drift and ensures
long-term stability in time-dependent simulations.
Third, Discrete Duality Finite Volume (DDFV) methods demonstrate
remarkable robustness on non-orthogonal, anisotropic, and severely
distorted meshes where traditional finite volume methods lose accuracy
or fail entirely~\cite{Domelevo-Omnes:2005,Boyer-Hubert:2008}.
Finally, dual grid discretizations achieve superior stability
properties through their inherent mathematical structure, with
rigorous analysis proving unconditional inf-sup stability for DDFV
methods~\cite{Boyer-Krell-Nabet:2015} and optimal convergence rates
for mimetic methods on general polyhedral
meshes~\cite{Brezzi-Lipnikov-Shashkov:2005}.

As a first step towards a dual grid approximation of more complex
application problems, such as those mentioned above, our approach for
the parabolic equation utilizes hyper-rectangular grids that match the
domain shape, maintaining two interconnected representations of the
solution.
One representation is defined at the vertices of the mesh
(vertex-centered) containing the primary solution values, the other at
the cell centers (cell-centered), providing a complementary
approximation that facilitates more accurate discretization of
differential operators.
These two representations work together during the discretization
process, where the information is transferred between them through
carefully designed operators that preserve the overall accuracy of the
scheme while creating a more robust numerical framework compared to
single-grid methods.
The discrete formulation of partial differential operations on this
dual grid representation transforms cell-based unknowns into
vertex-based unknowns (and vice versa), requiring specialized
formulations of gradient, divergence, and Laplacian operators that
maintain consistency across both representations and preserve
important mathematical properties of the continuous differential
operators.
This hybrid spatial discretization is finally coupled with first-order
explicit or implicit Euler schemes, or the second-order accurate
semi-implicit Crank-Nicolson method for time
integration~\cite{Marchuk:1990-book}.
Our methods ultimately yield an ODE system of the form:
\begin{align*}
  \frac{d\usV(\cdot;t)}{\dt} = \DeltaV\usV(\cdot;t) + \fs(\cdot;t),
\end{align*}
where $\DeltaV\usV(\cdot;t)$ is the vertex Laplace term built by
applying appropriate discretization formulas to $\usV(\cdot;t)$, while
$\fs(\cdot;t)$ is the source term evaluated at the mesh vertices at
time $t$.
The resulting schemes are second-order accurate in space and
first-order or second-order accurate in time, depending on the chosen
time-marching scheme.

\subsection{Tensor-train format and related work}
A conventional ``full-grid'' implementation of such schemes inevitably
faces the \emph{curse of dimensionality}, where memory and
computational costs grow exponentially with the number of
dimensions~\cite{Bellman:1961}.
Even for three-dimensional problems, calculations become prohibitively
expensive after a few mesh refinements.
Our primary contribution is reformulating this finite
difference/finite volume method into a more efficient algorithm using
the low-rank tensor-train format to mitigate the curse of
dimensionality.

The tensor-train (TT) format represents high-dimensional tensors as
products of lower-dimensional TT cores, each associated with a spatial
dimension~\cite{oseledets2011}.
For an $\NVx\times\NVy\times\NVz$ grid, the storage required by the TT
representation scales as $\calO\big(N\rs^2\big)$, where
$\Ns=\max(\NVx,\NVy,\NVz)$ and $\rs$ is the TT rank, compared to
$\calO\big(N^{3}\big)$ for conventional approaches. When $\rs\ll\Ns$,
this scaling yields substantial savings in memory and computational
cost.
It is important to distinguish between the tensor-train format
and the quantized tensor-train (QTT) format, which represents a
distinct extension of TT.
While TT operates directly on the original tensor dimensions, QTT
introduces additional (often non-physical) dimensions through binary
encoding of the original indices~\cite{khoromskij2011}.
This distinction is crucial, as most existing literature focuses on
either QTT-based approaches or applications to finite element methods
rather than finite difference/finite volume methods for single and
dual grid formulations.

\medskip
Previous work on tensor methods for partial differential equations
(PDEs) includes several distinct research directions.
Several studies~\cite{Kazeev2018,Kazeev2022,Markeeva2021} have
developed QTT-based finite element methods for elliptic problems, with
applications to multiscale diffusion and curvilinear domains.
These works primarily focus on steady-state equations rather than
time-dependent problems.
In Reference~\cite{Dolgov-Khoromskij-Oseledets:2012}, both TT and QTT
formats are applied to parabolic problems, particularly the
Fokker-Planck equation, using alternating linear schemes (ALS) as
their core numerical method.
Their approach differs fundamentally from ours, which employs
time-marching schemes directly on the TT cores.
Several papers have explored Krylov subspace methods for linear
systems with tensor product
structure~\cite{Ballani-Grasedyck:2012,Kressner-Tobler:2010,Oseledets-Dolgov:2012},
including TT-GMRES~\cite{Dolgov:2013} and general Krylov frameworks on
low-rank tensor
manifolds~\cite{Kressner-Steinlechner-Vandereycken:2016}.
While these works provide important algorithmic foundations, they do
not address the specific challenges of time-dependent parabolic
problems in the TT format using finite difference discretizations.
Various preconditioners have been developed for tensor-structured
problems~\cite{Bachmayr-Kazeev:2020,Dolgov-Khoromskij-Oseledets:2012}, though
most focus on general tensor formats or specific application domains
rather than time-dependent parabolic problems with hybrid spatial
discretizations.

\subsection{Novel contributions}
\label{subsec:contributions}

This work introduces several novel contributions that advance the
state-of-the-art in tensor-based numerical methods for time-dependent
PDEs.
All of these contributions are developed within the tensor-train
framework introduced by Oseledets~\cite{oseledets2011}.
In particular, our focus on finite difference/finite volume
discretizations for time-dependent parabolic problems in the TT format
presents unique challenges and opportunities that have not been
previously addressed in the literature.
Although we do not explore their potential application to the QTT
framework, such extensions are expected to be straightforward in most
cases and will be the focus of future research.
We summarize the key contributions of this work as follows.

\PARAGRAPH{$\mathbf{(i)}$ Novel hybrid spatial discretization in TT format.}
We present the first comprehensive framework that combines
cell-centered and vertex-centered finite difference approximations
within the tensor-train representation for three-dimensional parabolic
problems.
Unlike existing TT/QTT-based approaches, which typically focus on a
single discretization type, our dual-grid methodology enables more
flexible and accurate spatial approximations while maintaining the
computational efficiency of the TT format.
Unlike standard Laplace approximations represented by Kronecker sums
of the form $\matD^2 \otimes \matI \otimes \cdots \otimes \matI +
\cdots + \matI \otimes \cdots \otimes \matI \otimes \matD^2$, where
$\matD^2 = (1/h) \, \text{tridiag}(-1, 2, -1)$, our discretizations
involve transformations between cell-based and vertex-based unknowns
that require specialized reconstruction techniques and operators
working directly on TT cores.
This approach is related to the dual grid discretizations discussed
previously.
While we do not yet have a theoretical framework to analyze this
methodology in the tensor-train framework, our numerical experiments
demonstrate that the approach is \textit{rank-stable}, preserving the
rank of low-rank initial solutions through suitable rounding
procedures in the time-marching method.

\PARAGRAPH{$\mathbf{(ii)}$ Extension to variable meshes and domain remapping.}
We generalize our TT-based methodology to support non-uniform spatial
discretizations and remapped domains, expanding the applicability of
tensor-train methods beyond the uniform grid assumptions commonly
found in the literature.
This extension maintains second-order spatial accuracy while
preserving the computational benefits of the TT representation.

\PARAGRAPH{$\mathbf{(iii)}$ Direct implementation of multiple time-stepping schemes on TT cores.}
We redesign time-stepping algorithms—including explicit Euler,
implicit Euler, and Crank-Nicolson schemes—to operate directly on the
TT representation.
This approach ensures computational efficiency across the entire time
integration process by manipulating TT cores directly.
Unlike previous approaches that rely on alternating linear schemes
(ALS)~\cite{Dolgov-Khoromskij-Oseledets:2012}, our methods implement
classical time-marching schemes directly in the TT format, maintaining
the mathematical structure of the original schemes while exploiting
the computational advantages of the TT representation.

\PARAGRAPH{$\mathbf{(iv)}$ TT-specialized matrix-free preconditioned conjugate gradient method.}
To solve the linear systems arising from implicit and semi-implicit
time-stepping schemes, we reformulate the preconditioned conjugate
gradient (PCG) algorithm to operate directly on the TT cores of the
Krylov direction vectors.
To this end, we employ a \textit{matrix-free formulation} of the
matrix-vector product, so that the latter is provided by an
application of the discrete time operator to the Krylov directions.
Practically, we never build the time operator as a matrix in any given
full or tensor-compressed format.
This approach offers two key advantages:  

\vspace{-0.5\baselineskip}  
\begin{itemize}[itemsep=2.5pt, parsep=0pt]  
\item eliminates the need to convert Krylov direction vectors into
  full tensor representations during the iterative process;
\item avoids the explicit assembly or construction of the system
  matrix in both full and approximate TT representations,
  significantly reducing computational and memory overhead.
\end{itemize}  

\vspace{-0.5\baselineskip}  
\noindent  
To complete the PCG algorithm, we propose a specialized preconditioner
that acts directly on the TT cores by leveraging univariate finite
difference operators and exploiting the separable structure of the
discretization.
While previous works have explored preconditioners for
tensor-structured
problems~\cite{Bachmayr-Kazeev:2020,Dolgov-Khoromskij-Oseledets:2012},
our approach is specifically designed for the discrete operators
arising from our hybrid spatial discretization and focuses on
maintaining rank stability throughout the solution process.
By bypassing the explicit construction of the system matrix in TT
format, the method preserves the low-rank structure throughout the
preconditioning and solution process, while achieving optimal
computational scaling with refinement in space and time.

\PARAGRAPH{$\mathbf{(v)}$ Comprehensive algorithmic framework.}
We provide a complete algorithmic description, including a detailed
analysis of computational complexity, memory requirements, and
numerical accuracy.
In particular, we demonstrate that our TT reformulation preserves the
second-order accuracy of the hybrid finite difference approximation,
thereby ensuring theoretical guarantees for spatial and temporal
accuracy within the TT framework.
We also carefully describe technical implementation details, such as 
the setting of the boundary conditions, or the choice of threshold 
values, making it easy to implement our algorithms and replicate our 
results.

\medskip
In summary, the combination of these contributions results in
substantial practical improvements: memory requirements scale as
$\mathcal{O}(N r^2)$ instead of $\mathcal{O}(N^3)$ for traditional
methods, where $N$ is the grid size per dimension and $r$ is the TT
rank.
For problems where $r \ll N$, this represents
\emph{orders-of-magnitude improvements in both memory usage and
computational time}, making previously intractable three-dimensional
parabolic problems computationally feasible.
This methodology opens new possibilities for solving large-scale
multidimensional problems in the area of finite difference/finite
volume methods that were previously computationally prohibitive with
conventional approaches.

\subsection{Paper's outline}
The remainder of this paper is organized as follows:
Section~\ref{sec:TTformat} presents the fundamental mathematical
framework, introducing the notation and essential concepts of the
tensor-train format that underpins our methodology.
Section~\ref{sec:FVM-regular} develops the finite difference/finite
volume discretization scheme and its reformulation within the
tensor-train framework for the numerical approximation of the system
defined by Equation~\eqref{eq:problem:strong}.
In Section~\ref{sec:FVM-remapped}, we generalize our approach by
extending the discretization method to handle problems on variable
mesh sizes and remapped domains, maintaining the framework established
in Section~\ref{sec:FVM-regular}.
Section~\ref{sec:time_integration} details the implementation of the
different time marching techniques, including the preconditioned
conjugate gradient method for solving the linear systems that arise
from implicit and semi-implicit time discretizations.
Section~\ref{sec:numerical} demonstrates the effectiveness and
efficiency of our proposed methodology through comprehensive numerical
experiments and performance analysis.
Section~\ref{sec:conclusions} summarizes our findings, discusses their
implications, and outlines potential directions for future research.


\section{Tensor-train format: notation, basic definitions, background properties}
\label{sec:TTformat}

We shorty recall here some notation, basic definitions, and a few
technicalities about tensors and the tensor-train format, referring to the review article~\cite{Kolda-Bader:2009} and the pioneering
papers~\cite{oseledets2010a,oseledets2011} for a
detailed exposition.
A \emph{$\DIM$-dimensional tensor} $\tA$ is a multi-dimensional array
with $\DIM$ indices, e.g.,
$\tA=\SET{\tA(\iss{1},\iss{2},\ldots,\iss{\DIM})},
\iss{\ell}=1,2,\ldots,\nss{\ell}, \ell=1,2,\ldots,\DIM$.
Since we are interested in solving three-dimensional (3D) problems, we
restrict all basic definitions to $0\leq\DIM\leq3$, and use the Latin
indices $i,j,k$ instead of $\iss{1},\iss{2},\iss{3}$.
Technically, scalar quantities, vectors, and matrices are special
cases of $\DIM$-way tensors for $\DIM=0,1,2$, respectively, but
throughout the paper, we prefer to use a distinct notation for them
for clarity of exposition.
We denote a \emph{scalar quantity} using normal, lower-case fonts,
e.g. $a,b,c$, etc; a \emph{vector field} using bold, lower-case fonts,
e.g., $\av,\bv,\ccv$, etc; a \emph{matrix} using normal, upper case
fonts, e.g, $\As,\Bs,\Cs$, etc; a \emph{$3$-dimensional tensor} using
calligraphic, upper case fonts, e.g., $\tA,\tB,\tC$, etc.
We address the entries of $\DIM$-dimensional quantities for
$\DIM\geq1$ by using a MATLAB-like notation; accordingly, the symbol
``:'' denotes a \emph{free index}.
So, $\as(i)$, $\As(i,j)$, $\tA(i,j,k)$ are the $i$-th component of
vector $\av$, the $(i,j)$-th component of matrix $\As$, and the
$(i,j,k)$-th component of tensor $\tA$, respectively.

It is convenient to introduce the concepts of \emph{fibers} and
\emph{slices}.
A fiber is the vector that we obtain from a tensor by fixing all
indices but one; a slice is the matrix that we obtain by fixing all
indices but two.
Therefore, a 3D tensor $\tA(i,j,k)$ has three distinct fibers, e.g.,
$\tA(:,j,k)$, $\tA(i,:,k)$, and $\tA(i,j,:)$, and three distinct
slices, e.g., $\tA(i,:,:)$, $\tA(:,j,:)$, and $\tA(:,:,k)$.

\medskip
We say that a 3D tensor
$\tA\in\REAL^{\nss{1}\times\nss{2}\times\nss{3}}$ is in
\emph{tensor-train format} if there exist three 3D tensors, called
\emph{tensor-train cores} and denoted by
$\tAss{\ell}\in\REAL^{\rss{\ell-1}\times\nss{\ell}\times\rss{\ell}}$,
for $\ell=1,2,3$ and with $\rss{0}=\rss{3}=1$, such that
\begin{align}
  \tA(i,j,k)
  = \sum_{\alpha_{1}=1}^{\rss{1}}\sum_{\alpha_{2}=1}^{\rss{2}}
  \tAss{1}(          i,\alpha_{1})
  \tAss{2}(\alpha_{1},j,\alpha_{2})
  \tAss{3}(\alpha_{2},k).
  \label{eq:TT:def}
\end{align}
Since $\rss{0}=\rss{3}=1$, the two cores $\tAss{1}$ and $\tAss{3}$ are
indeed matrices, so we write them with two indices
in~\eqref{eq:TT:def}, but still keeping the tensor notation.
We can reformulate definition~\eqref{eq:TT:def} by introducing three
slice matrices
$\Ass{1}(i)=\tAss{1}(:,i,:)\in\REAL^{\rss{0}\times\rss{1}}$,
$\Ass{2}(j)=\tAss{2}(:,j,:)\in\REAL^{\rss{1}\times\rss{2}}$, and
$\Ass{3}(k)=\tAss{3}(:,k,:)\in\REAL^{\rss{2}\times\rss{3}}$,
which are parametrized with the spatial indices $i$, $j$, and $k$,
such that
\begin{align*}
  \tA(i,j,k) = \Ass{1}(i)\Ass{2}(j)\Ass{3}(k),
\end{align*}
for every possible combination of the indices $i$, $j$, $k$.
The integers $\rss{1},\rss{2}$ represent the \textit{TT ranks} of
$\tA{}$, that is, the internal sizes of tensors
$\tAss{1},\tAss{2},\tAss{3}$.
It is also convenient to define an upper bound of these ranks as
$\rs=\max\left\{\rss{\ell}\right\}_{\ell}$.

In the next section, we will introduce two distinct types of 3D
tensors, which are logically associated with the vertices and the
cells of a 3D mesh, respectively.
To emphasize the nature of such tensors, instead of $(i,j,k)$ we will
use the index triple $(\iV,\jV,\kV)$ in the first case, and the index
triple $(\ic,\jc,\kc)$ in the second.
This double index notation of vertices and cells avoids using
``half-indices''; usually, we might denote two consecutive vertices
as, for example, ``$i$'' and ``$i+1$'' along the $\dirX$ direction and
the cell between them as ``$i+\frac12$''.
However, this notation is very cumbersome when we need to index the
cores of the tensor-train representation of a cell grid function as these indices should be integers and not half integers.
For this reason, we prefer using this double-index notation, which is
completely equivalent to the more familiar one once we observe that the $\iV$-th vertex and the $\ic$-th cell correspond to the $i$-th
vertex and the $(i+\frac12)$-th cell, respectively, if we assume that
$i=\iV=\ic$, as it is discussed in the next section.

A very efficient tensor-train decomposition algorithm called \TTSVD{}
is available from Ref.~\cite{oseledets2011}, which is based on a
sequence of Singular Value Decomposition (SVD) of auxiliary matrices.
The most crucial properties of this algorithm are those related to the TT ranks and the ranks of the so-called \textit{matrix unfoldings}.
A three-dimensional tensor has two distinct matrix unfoldings.
The first unfolding is the matrix
$\As_1\in\REAL^{\nss{1}\times(\nss{2}\nss{3})}$ whose row index
corresponds to the first tensor index $i$, and whose column index is a
bijective remap, say $\eta_1$, of the tensor indices $(j,k)$, so that
$\As_1(i,\eta_1(j,k))=\tA(i,j,k)$.
The second unfolding is the matrix
$\As_2\in\REAL^{(\nss{1}\nss{2})\times\nss{3}}$ whose row index is a
bijective remap, say again $\eta_2$, of the tensor indices $(i,j)$,
and column index is equal to $\kV$ so that
$\As_2(\eta_2(i,j),k)=\tA(i,j,k)$.

According to~\cite[Theorem~2.1]{oseledets2011}, \TTSVD{} can compute
an exact representation of the tensor $\tA$ in the TT format with
TT ranks $\rss{\ell}$ satisfying
$\rss{\ell}=\operatorname{rank}(\As_{\ell})$ for both $\ell=1,2$.
According to~\cite[Corollary~2.4]{oseledets2011}, if we select smaller
ranks $r_{\ell}'\leq\operatorname{rank}(\As_{\ell})$, then the best
approximation
\begin{align*}
 \tATT_{\text{best}} 
 = \arg\inf_{\,\tATT\in\cTT(\rv'_\ell)\,}\|\tA - \tATT\|_{\tF}
\end{align*}
\vspace{-0.65\baselineskip}

\noindent
in the TT format of $\tA$
exists with such $r_{\ell}'$ as internal ranks.
Here, the \emph{infimum} is taken over the set of tensors
$\cTT(\rv_\ell')$ in tensor-train format having internal ranks
$\rv_\ell'=(\rss{1}',\rss{2}')$.
The \TTSVD{} algorithm balances computational efficiency and
approximation quality since it returns a \emph{quasi-optimal best
approximation} $\tATT$ in the sense that
\begin{align*}
  \| \tA - \tATT \|_{\tF} \leq 
  \Cs_{\TT}(\DIM)
  \|\tA - \tATT_{\text{best}}\|_{\tF}.
\end{align*}
For $\DIM=3$, the inequality constant is $\Cs_{\TT}=\sqrt{2}$ since the general expression is
$\Cs_{\TT}(\DIM)=\sqrt{\DIM-1}$~\cite{oseledets2010a}.

Multi-linear algebra operations, such as addition, multiplication by a
scalar, element-wise (Hadamard) multiplication, contraction with
matrices and vectors are possible through a straightforward implementation that facilitates tensor
calculations~\cite{oseledets2011}.
However, these operations may increase the TT ranks, and an
application of the so-called \textit{rounding}
algorithm~\cite{oseledets2011} is required to control the rank growth.
We denote this operation with the symbol $\RNDG(\,\cdot\,)$.
Included in the $\RNDG$ procedure, there is a threshold $\toll$ which
ensures that the ``rounded'' tensor-train will approximate the original one accurately.
The choice of $\toll$ is crucial for the good behavior of the method.
Smaller $\toll$ values lead to higher TT ranks, resulting in more
accurate but more expensive computations, and vice versa.

A sampling algorithm based on the \MAXVOL{} row and column selection
strategy~\cite{Goreinov-Oseledets-etal:2010} is available for
constructing an approximate representation of a tensor in the TT
format without precomputing and storing the tensor, see,
e.g.,~\cite{oseledets2010a}.
We use this algorithm to ensure efficient processing to the treatment
of the right-hand side term $\fs$.



\section{Tensor-train reformulation of finite difference operators on 3D cartesian grids}
\label{sec:FVM-regular}

In this section, we focus exclusively on the spatial discretization.
Therefore, for clarity and simplicity, we will omit the time
dependence of spatially varying fields, i.e., we will discretize
$\us(\xv)$ and its derivatives in space instead of $\us(\xv,\ts)$.

\subsection{Regular equispaced meshes and grid functions}
\label{subsec:FVM-regular:mesh}
We consider the computational domain
$\Omega=(\ass{x},\bss{x})\times(\ass{y},\bss{y})\times(\ass{z},\bss{z})$,
for the three pairs of real values $\ass{\ell}$ and $\bss{\ell}$, with
$-\infty<\ass{\ell}<\bss{\ell}<\infty$, $\ell\in\{x,y,z\}$.
We partition $\overline{\Omega}$ along the directions $\dirX$, $\dirY$
and $\dirZ$
with constant mesh size steps $\hx$, $\hy$, $\hz$, respectively, and
set the ``characteristic mesh size'' $\hh=\min(\hx,\hy,\hz)$.
This partitioning introduces $\NVx$, $\NVy$, and $\NVz$ vertices and
$\NCx=\NVx-1$, $\NCy=\NVy-1$, $\NCz=\NVz-1$ univariate closed
intervals along $\dirX$, $\dirY$ and $\dirZ$, respectively.
The mesh size steps satisfy the conditions:
$\NCx\hx=\bss{x}-\ass{x}$, $\NCy\hy=\bss{y}-\ass{y}$, and
$\NCz\hz=\bss{z}-\ass{z}$.

We label each mesh vertex with the index triple $(\iV,\jV,\kV)$, so
that the $(\iV,\jV,\kV)$-th vertex, also denoted as $\V(\iV,\jV,\kV)$,
has coordinates $\xvV=\big(\xsV(\iV),\ysV(\jV),\zsV(\kV)\big)$, where
\begin{align*}
  \xsV(\iV)=\ass{x}+\iV\hx \quad\textrm{for~}\iV=0,1,\ldots,\NVx-1,\\
  \ysV(\jV)=\ass{y}+\jV\hy \quad\textrm{for~}\jV=0,1,\ldots,\NVy-1,\\
  \zsV(\kV)=\ass{z}+\kV\hz \quad\textrm{for~}\kV=0,1,\ldots,\NVz-1.  
\end{align*}

We label each mesh cell with the index triple $(\ic,\jc,\kc)$, so that
the $(\ic,\jc,\kc)$-th cell, also denoted as $\C(\ic,\jc,\kc)$, is
given by the tensor product of the univariate closed intervals along
$\dirX$, $\dirY$ and $\dirZ$
\begin{align*}
  \C(\ic,\jc,\kc) &=
  \big[\xsV(\iV,\jV,\kV), \xsV(\iV+1,\jV,\kV)\big] \times
  \big[\ysV(\iV,\jV,\kV), \ysV(\iV,\jV+1,\kV)\big] \times
  \\[0.5em]
  &\hspace{7cm}\times
  \big[\zsV(\iV,\jV,\kV), \zsV(\iV,\jV,\kV+1)\big].
\end{align*}
Cell indices run as $\ic=0,1,\ldots\NCx-1$, $\jc=0,1,\ldots\NCy-1$,
$\kc=0,1,\ldots\NCz-1$.
Cell $ \C(\ic,\jc,\kc)$ can also be defined as the convex envelope of
its eight vertices:
\begin{align}
  \C(\ic,\jc,\kc)
  &= \textrm{convex~envelope~of}\nonumber\\[0.5em]
  &\,\Big\{
  \V(\iV,\jV,\kV), \V(\iV,\jV+1,\kV), \V(\iV,\jV+1,\kV+1\big),
  \V(\iV,\jV,\kV+1\big),\nonumber\\
  &\qquad
  \V(\iV+1,\jV,\kV), \V(\iV+1,\jV+1,\kV), \V(\iV+1,\jV+1,\kV+1), \V(\iV+1,\jV,\kV+1)
  \Big\}.\label{eq:cell:convex:def}
\end{align}
Consequently, we can identify cell $\C(\ic,\jc,\kc)$ as the cell whose
``first'' vertex in definition~\eqref{eq:cell:convex:def} is
$\V(\iV,\jV,\kV)$, provided that $\iV=\ic$, $\jV=\jc$, and $\kV=\kc$.
This fact introduces a precise and well-defined bijective
correspondence between the vertex and the cell numbering systems.
The cell center has coordinates $\xvC=(\xsC(\ic),\ysC(\jc),\zsC(\kc))$,
where
\begin{align*}
  \xsC(\ic)=\ass{x}+\left(\ic+\frac12\right)\hx=\frac{\xs(\iV)+\xs(\iV+1)}{2},\\
  \ysC(\jc)=\ass{y}+\left(\jc+\frac12\right)\hy=\frac{\ys(\jV)+\ys(\jV+1)}{2},\\
  \zsC(\kc)=\ass{z}+\left(\kc+\frac12\right)\hz=\frac{\zs(\kV)+\zs(\kV+1)}{2}.
\end{align*}

We represent a scalar field $\us(\xs,\ys,\zs)$ as a \emph{vertex grid
function} $\usV=\SET{\usV(\iV,\jV,\kV)}$ or a \emph{cell grid
function} $\usC=\SET{\usC(\ic,\jc,\kc)}$.
These grid functions are three-dimensional tensors
regardless of being associated with vertices or cells.
Depending on the context, we may interpret each value
$\usV(\iV,\jV,\kV)$ as either $\us(\xv_{\V})$, i.e., the value of
$\us(\xv)$ sampled at the $(\iV,\jV,\kV)$-th vertex with coordinates
$\xvV$, or as an approximation of it.
Similarly, we may interpret each value $\usC(\ic,\jc,\kc)$ as either
$\us(\xv_{\C})$, i.e., the value of $\us(\xv)$ sampled at the center
of the $(\ic,\jc,\kc)$-th cell with coordinates $\xvC$, or as an
approximation of it.
A second interpretation is also possible where this quantity is either the cell average or an approximation of the cell average of
$\us(\xv)$ over cell $\C(\ic,\jc,\kc)$.
The first interpretation is related to a finite difference setting, the second one to a finite volume setting, but on orthogonal,
cartesian grids, they are equivalent up to an error scaling
like $\calO(\hh^2)$, and may eventually lead to the same approximation
formulas.
Throughout this paper, we leave the interpretation of
$\usC(\ic,\jc,\kc)$ unspecified, thus referring implicitly to both, although sometimes it may look like we have some preference for the finite
difference one as we often implicitly refer to the concept of a
``pointwise stencil''.

\subsection{First derivatives of vertex grid functions}
\label{subsec:FVM-regular:first_derivative}
In our setting, the discrete derivative of a vertex grid function is a
cell grid function.
Let $\usVx$, $\usVy$, and $\usVz$ denote the three components of the discrete gradient of the grid function $\usV$.
They are associated with cell $\C(\ic,\jc,\kc)$ and defined as
follows:
\begin{align}
  \usVx(\ic,\jc,\kc)
  &:= \frac{1}{4\hx}\sum_{n=0}^{1}\sum_{p=0}^{1}\Big(\usV(\iV+1,\jV+n,\kV+p) - \usV(\iV,\jV+n,\kV+p)\Big),
  \label{eq:first:X}\\[0.5em]
  \usVy(\ic,\jc,\kc)
  &:= \frac{1}{4\hy}\sum_{m=0}^{1}\sum_{p=0}^{1}\Big(\usV(\iV+m,\jV+1,\kV+p) - \usV(\iV+m,\jV,\kV+p)\Big),
  \label{eq:first:Y}\\[0.5em]
  \usVz(\ic,\jc,\kc)
  &:= \frac{1}{4\hz}\sum_{m=0}^{1}\sum_{n=0}^{1}\Big(\usV(\iV+m,\jV+n,\kV+1) - \usV(\iV+m,\jV+n,\kV)\Big).
  \label{eq:first:Z}
\end{align}

\begin{figure}[htbp]
    \centering
    \includegraphics[width=0.6\linewidth]{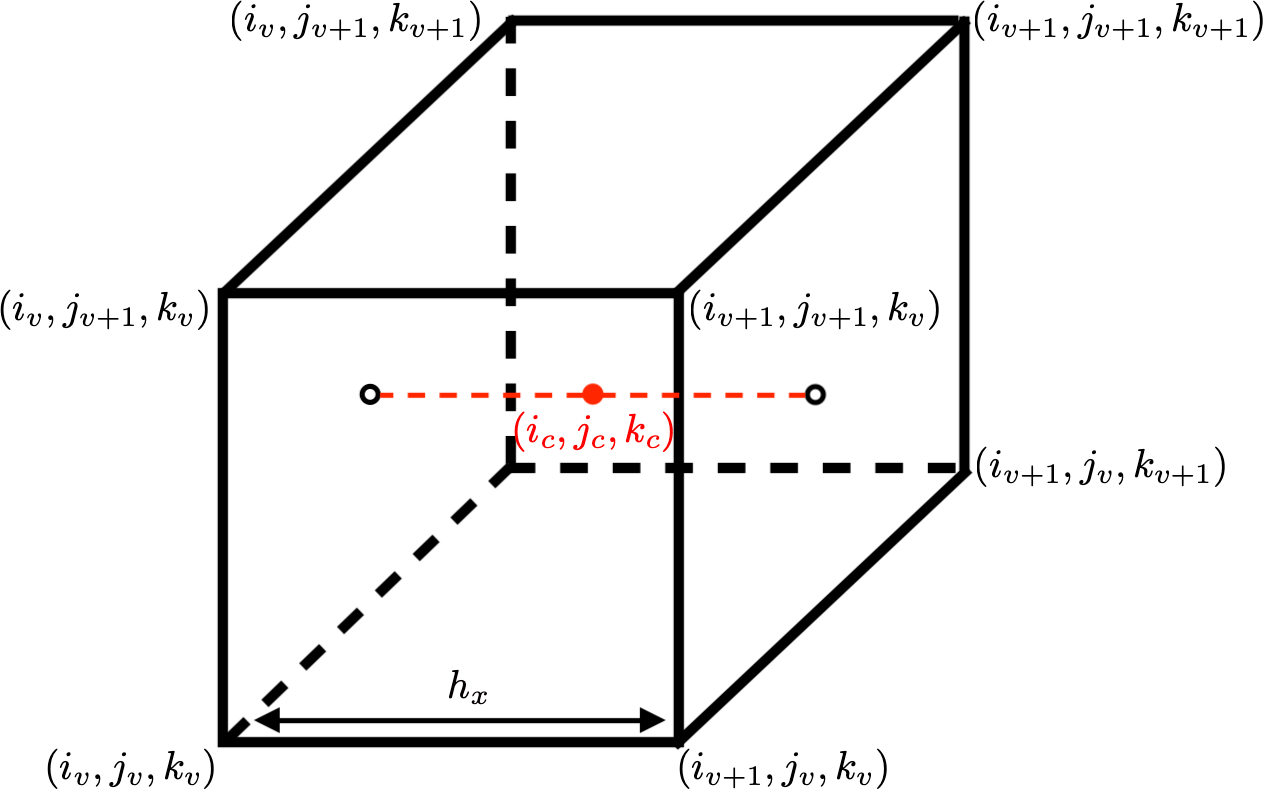}
    \caption{Computation of $\usVx(\ic,\jc,\kc)$.}
    \label{fig:stencil_dx}
\end{figure}

Note that the difference operations are applied to the vertex indices, but the result is associated with the cell identified by taking
$\ic=\iV$, $\jc=\jV$, $\kc=\kV$ according to the numbering system of
cells and vertices introduced in
Section~\ref{subsec:FVM-regular:mesh}, see
Figure~\ref{fig:stencil_dx}.
The association to a cell is reflected by the superscript ``$c$'' in
the discrete difference operators $\partial^{c}_{x}$,
$\partial^{c}_{y}$, and $\partial^{c}_{z}$.
Then, assuming that
$\usV(\iV+m,\jV+n,\kV+p)=\us(\xsV(\iV+m),\ysV(\jV+n),\zsV(\kV+p))$
with $m,n,p=0,1$, a straightforward Taylor expansion implies that
\begin{equation}
  \label{eq:2nd-order:accuracy}
  \begin{aligned}
    &
    \left(\frac{\partial\us}{\partial\xs}\right)_{\C}=\usVx(\ic,\jc,\kc)+\calO(\hh^2),\qquad
    \left(\frac{\partial\us}{\partial\ys}\right)_{\C}=\usVy(\ic,\jc,\kc)+\calO(\hh^2),\\[1.em]
    &
    \left(\frac{\partial\us}{\partial\zs}\right)_{\C}=\usVz(\ic,\jc,\kc)+\calO(\hh^2),
  \end{aligned}
\end{equation}
where, to ease the notation, the symbols
$\big({\partial\us}\slash{\partial\xs}\big)_{\C}$,
$\big({\partial\us}\slash{\partial\ys}\big)_{\C}$, and
$\big({\partial\us}\slash{\partial\zs}\big)_{\C}$ denote the
corresponding first derivatives evaluated at the cell center
$\big(\xsC(\ic),\ysC(\jc),\zsC(\kc)\big)$.
We report the details of this calculation in~\ref{sec:appx:Taylor_expansion_1st_derivatives}.

\medskip
Formulas~\eqref{eq:first:X}, \eqref{eq:first:Y},
and~\eqref{eq:first:Z} are rather standard in the finite difference
and finite volume literature. 
However, they suffer the curse of dimensionality as the
computational complexity of computing $\usVx$, $\usVy$, and $\usVz$ on
a 3D, $(\NVx\times\NVy\times\NVz)$-sized grid scales like
$\calO(\Ns^3)$ with $\Ns=\max(\NVx,\NVy,\NVz)$ as already noted in the
introduction.
This ``bad'' scaling makes calculations expensive and almost
unfeasible even for relatively small values of $\Ns$, e.g.,
between $10^2$ and $10^3$.
Herein, we are interested in reformulating such first derivatives in
the tensor-train format.
Let
\begin{align*}
  \usVTT(\iV,\jV,\kV)
  = \sum_{\alpha_1=1}^{\rss{1}}\sum_{\alpha_2=1}^{\rss{2}} \tUss{1}(\iV,\alpha_1)\tUss{2}(\alpha_1,\jV,\alpha_2)\tUss{3}(\alpha_2,\kV)
  = \Uss{1}(\iV)\Uss{2}(\jV)\Uss{3}(\kV)
\end{align*}
be the TT representation of the grid function $\usV$, using cores
$\tUss{\ell}$ and matrix slices $\Uss{\ell}$ for $\ell=1,2,3$
according to the notation and definitions discussed in
Section~\ref{sec:TTformat}.
At the moment, we do not need to specify if this is an ``exact'' or an
``approximate'' representation, nor if we obtained the grid function
$\usV$ (represented as $\usVTT$) by sampling some smooth function
$\us(\xv)$ at the vertices of a given grid, or if $\usV$ is the result
of some complex calculations.
For convenience of exposition, we assume that this is an exact
representation so that it holds that
$\usV(\iV,\jV,\kV)=\usVTT(\iV,\jV,\kV)$ for every meaningful index
triple $(\iV,\jV,\kV)$.
To compute the first derivative of $\usV$ along the direction $\dirX$
in the TT format, i.e., the cell grid function $\usVTTx(\ic,\jc,\kc)$,
we apply definition~\eqref{eq:first:X}, and we find that
\begin{align}
  \usVTTx(\ic,\jc,\kc)
  &:= \frac{1}{4\hx}\sum_{n=0}^{1}\sum_{p=0}^{1}\Big(\usVTT(\iV+1,\jV+n,\kV+p) - \usVTT(\iV,\jV+n,\kV+p)\Big).
  \label{eq:first:X:TT}
\end{align}
Then, we expand the first summation on the right
of~\eqref{eq:first:X:TT}, and using the tensor-train definition, we
rewrite it as
\begin{align*}
  &\sum_{n=0}^{1}\sum_{p=0}^{1}\usVTT(\iV+1,\jV+n,\kV+p)
  =
  \usVTT(\iV+1,\jV,  \kV  ) +
  \usVTT(\iV+1,\jV+1,\kV  ) \\[-0.75em]
  & \hspace{6cm}
  + \usVTT(\iV+1,\jV+1,\kV+1)
  + \usVTT(\iV+1,\jV,  \kV+1)
  \\[1em] &\qquad
  =
  \Uss{1}(\iV+1)\Uss{2}(\jV  )\Uss{3}(\kV  ) +
  \Uss{1}(\iV+1)\Uss{2}(\jV+1)\Uss{3}(\kV  ) +
  \Uss{1}(\iV+1)\Uss{2}(\jV+1)\Uss{3}(\kV+1)
  \\[0.5em] &\hspace{6cm}
  +\Uss{1}(\iV+1)\Uss{2}(\jV  )\Uss{3}(\kV+1)
  \\[0.75em] &\qquad
  =
  \Uss{1}(\iV+1)\,\Big(\Uss{2}(\jV)+\Uss{2}(\jV+1)\Big)\,\Big(\Uss{3}(\kV) +\Uss{3}(\kV+1)\Big).
\end{align*}
We repeat the same calculation for the second summation term on the
right of Eq.~\eqref{eq:first:X:TT}, we take the difference between
them and substitute the result into Eq.~\eqref{eq:first:X:TT}.
We find that
\begin{align*}
  \usVTTx(\ic,\jc,\kc)
  &= \Usz{1}{x}(\ic)\,\Usz{2}{x}(\jc)\,\Usz{3}{x}(\kc)\\
  &:=
  \Bigg(\frac{\Uss{1}(\iV+1)-\Uss{1}(\iV  )}{\hx}\Bigg)\,
  \Bigg(\frac{\Uss{2}(\jV  )+\Uss{2}(\jV+1)}{2}  \Bigg)\,
  \Bigg(\frac{\Uss{3}(\kV  )+\Uss{3}(\kV+1)}{2}  \Bigg),  
\end{align*}
where we identify the vertex index pairs $(\iV,\iV+1)$, $(\jV,\jV+1)$,
$(\kV,\kV+1)$, with the cell indices $\ic,\jc,\kc$, and the cores
$\Usz{1}{x}$, $\Usz{2}{x}$ and $\Usz{3}{x}$ parametrized with
$(\ic,\jc,\kc)$ with the cores of the first derivative of $\usVTT$
along the direction $\dirX$.
The superscript ``$x$'' in the core symbols $\Usz{\ell}{x}$,
$\ell=1,2,3$ refers to the direction of the derivative.
We can similarly write the other two derivatives:
\begin{align*}
  \usVTTy(\ic,\jc,\kc)
  &= \Usz{1}{y}(\ic)\,\Usz{2}{y}(\jc)\,\Usz{3}{y}(\kc)\\
  &:=
  \Bigg(\frac{\Uss{1}(\iV  )+\Uss{1}(\iV+1)}{2}  \Bigg)\,
  \Bigg(\frac{\Uss{2}(\jV+1)-\Uss{2}(\jV  )}{\hy}\Bigg)\,
  \Bigg(\frac{\Uss{3}(\kV  )+\Uss{3}(\kV+1)}{2}  \Bigg),  
\end{align*}
and
\begin{align*}
  \usVTTz(\ic,\jc,\kc)
  &= \Usz{1}{z}(\ic)\,\Usz{2}{z}(\jc)\,\Usz{3}{z}(\kc)\\
  &:=
  \Bigg(\frac{\Uss{1}(\iV  )+\Uss{1}(\iV+1)}{2}  \Bigg)\,
  \Bigg(\frac{\Uss{2}(\jV  )+\Uss{2}(\jV+1)}{2}  \Bigg)\,
  \Bigg(\frac{\Uss{3}(\kV+1)-\Uss{3}(\kV  )}{\hz}\Bigg).
\end{align*}
A crucial point of this contruction is that no additional errors are
introduced in the TT representation of the derivatives other than
those of the finite difference approximation as
in~\eqref{eq:2nd-order:accuracy}.

\subsection{Second derivatives of vertex grid functions}
\label{subsec:FVM-regular:second_derivative}
\begin{figure}[htbp]
    \centering
    \includegraphics[width=0.8\linewidth]{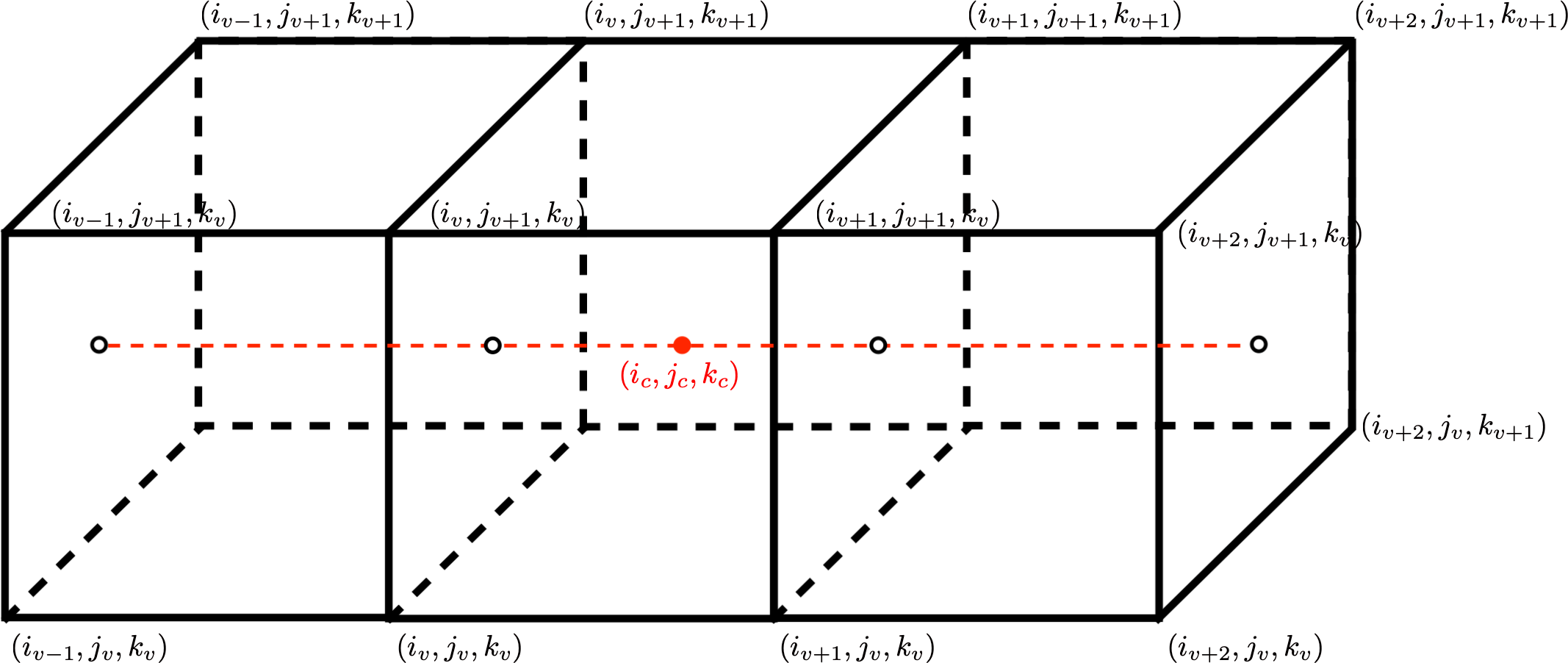}
    \caption{Computation of $\usVxx(\ic,\jc,\kc)$.}
    \label{fig:stencil_dxx}
\end{figure}

To define the discrete second derivatives of $\usVTT$, the vertex grid
function $\usV$ in the TT format representation, we proceed by still
assuming that the TT representation is exact, i.e.,
$\usVTT(\iV,\jV,\kV)=\usV(\iV,\jV,\kV)$, and using a similar argument
as for the first derivatives.
So, we first start from the discrete second derivative formulas for
$\usV$ at cell $\C(\ic,\jc,\kc)$ along $\dirX$, $\dirY$ and $\dirZ$,
which are given by 
\begin{align}
  \usVxx(\ic,\jc,\kc)  = \frac{1}{8\hx^2}
  \sum_{m=0}^{3}\eta_{m}\,\sum_{n=0}^{1}\sum_{p=0}^{1} \xi_{n}\xi_{p} \usV(\iV-1+m,\jV+n,\kV+p),
  \label{eq:second:X}\\[0.5em]
  \usVyy(\ic,\jc,\kc)  = \frac{1}{8\hy^2}
  \sum_{m=0}^{3}\eta_{m}\,\sum_{n=0}^{1}\sum_{p=0}^{1} \xi_{n}\xi_{p} \usV(\iV+n,\jV-1+m,\kV+p),
  \label{eq:second:Y}\\[0.5em]
  \usVzz(\ic,\jc,\kc)  = \frac{1}{8\hz^2}
  \sum_{m=0}^{3}\eta_{m}\,\sum_{n=0}^{1} \sum_{p=0}^{1} \xi_{n}\xi_{p} \usV(\iV+n,\jV+p,\kV-1+m),
  \label{eq:second:Z}
\end{align}
where $\xi_{0}=-1$, $\xi_{1}=1$, $\eta_{0}=\eta_{3}=1$ ,
$\eta_{1}=\eta_{2}=-1$, see Figure~\ref{fig:stencil_dxx}.

Assuming again that at every mesh vertex $\V(\iV,\jV,\kV)$ it holds
that $\usV(\iV,\jV,\kV)=\us(\xsV(\iV),\ysV(\jV),\zsV(\kV))$, a
straightforward Taylor expansion implies that
\begin{equation}
    \label{eq:2nd-order:accuracy:2nd-derivatives}
    \begin{aligned}
      &
      \left(\frac{\partial^{2}\us}{\partial\xs^2}\right)_{\C} =
      \usVxx(\iV,\jV,\kV) + \mathcal{O}(\hh^2),\quad
      \left(\frac{\partial^{2}\us}{\partial\ys^2}\right)_{\C} =
      \usVyy(\iV,\jV,\kV) + \mathcal{O}(\hh^2),\\[0.5em]
      &
      \left(\frac{\partial^{2}\us}{\partial\zs^2}\right)_{\C} =
      \usVzz(\iV,\jV,\kV) + \mathcal{O}(\hh^2),
    \end{aligned}
\end{equation}
where again $(\partial^2\us\slash{\partial\xs^2})_{\C}$,
$(\partial^2\us\slash{\partial\ys^2})_{\C}$, and
$(\partial^2\us\slash{\partial\zs^2})_{\C}$ are the second derivatives
of $\us$ computed at the center of cell $\C(\ic,\jc,\kc)$.
As reflected by the term $\mathcal{O}(\hh^2)$ these formulas provide a
second-order accurate approximation of the second derivatives of $\us$
along $\dirX$, $\dirY$ and $\dirZ$.
More details 
are reported in~\ref{sec:appx:Taylor_expansion_2nd_derivatives}.

Consider the discrete derivative along $\dirX$ given
by~\eqref{eq:second:X} (the same argument applies to the other
derivatives).
We can apply this formula directly to $\usVTT(\iV,\jV,\kV)$ since we
assume that the TT format representation of $\usV$ is exact, i.e.,
$\usVTT(\iV,\jV,\kV)=\usV(\iV,\jV,\kV)$.
As for the first derivatives, we must now see what this formula
implies for the cores of $\usVTT$.
For $\iV-1+m$ with $m=0,1,2,3$, by using the values of coefficients
$\xi_n,\xi_p=\pm1$, and expanding the double internal summation we
obtain:
\begin{align*}
  &\sum_{n=0}^{1}\sum_{p=0}^{1} \xi_{n}\xi_{p} \usVTT(\iV-1+m,\jV+n,\kV+p)\\[0.5em]
  &\,\,=\usVTT(\iV-1+m, \jV, \kV) - \usVTT(\iV-1+m, \jV+1, \kV) - \usVTT(\iV-1+m, \jV, \kV+1) \\
  &\hspace{9.5cm}+ \usVTT(\iV-1+m, \jV+1, \kV+1)\\[1em]
  &\,\,=
   \Uss{1}(\iV-1+m)\Uss{2}(\jV  )\Uss{3}(\kV)
  -\Uss{1}(\iV-1+m)\Uss{2}(\jV+1)\Uss{3}(\kV  )\\[0.25em]
  &\qquad
  -\Uss{1}(\iV-1+m)\Uss{2}(\jV  )\Uss{3}(\kV+1)
  +\Uss{1}(\iV-1+m)\Uss{2}(\jV+1)\Uss{3}(\kV+1)\\[0.5em]
  &\,\,=
  \Uss{1}(\iV-1+m)\,\Big(\Uss{2}(\jV)+\Uss{2}(\jV+1)\Big)\,\Big(\Uss{3}(\kV)+\Uss{3}(\kV+1)\Big).
\end{align*}
Then, expanding the first summation term in~\eqref{eq:second:X}, using
the values of coefficient $\eta_m$ for $m=0,1,2,3$, and reordering the matrix
slices $\Uss{1}(\iV-1+m)$ for various $m$, we obtain:
\begin{align*}
  &\sum_{m=0}^{3}\eta_{m}\,\sum_{n=0}^{1}\sum_{p=0}^{1} \xi_{n}\xi_{p} \usVTT(\iV-1+m,\jV+n,\kV+p)\\
  &=\Bigg[\sum_{m=0}^{3}\eta_{m}\Uss{1}(\iV-1+m)\Bigg]\,\Big(\Uss{2}(\jV)+\Uss{2}(\jV+1)\Big)\,\Big(\Uss{3}(\kV)+\Uss{3}(\kV+1)\Big)\\
  &=\Big[ \Big(\Uss{1}(\iV+2)-\Uss{1}(\iV+1)\Big) - \Big(\Uss{1}(\iV)-\Uss{1}(\iV-1)\Big) \Big]
  \Big(\Uss{2}(\jV)+\Uss{2}(\jV+1)\Big)
  \Big(\Uss{3}(\kV)+\Uss{3}(\kV+1)\Big).
\end{align*}
Finally, we obtain the expression of the second discrete derivative of
$\usVTT$, which has the form:
\begin{align*}
  \usVTTxx(\ic,\jc,\kc)
  &= \Usz{1}{xx}(\ic)\,\Usz{2}{xx}(\jc)\,\Usz{3}{xx}(\kc)\\[0.5em]
  &:=  
  \frac{1}{2\hx}
  \left( \frac{\Uss{1}(\iV+2)-\Uss{1}(\iV+1)}{\hx} - \frac{\Uss{1}(\iV)-\Uss{1}(\iV-1)}{\hx} \right)
  \left( \frac{\Uss{2}(\jV)+\Uss{2}(\jV+1)}{2} \right)\times\\
  &\hspace{9cm}
  \times\left( \frac{\Uss{3}(\kV)+\Uss{3}(\kV+1)}{2} \right),
\end{align*}
where, as for the first derivative, we associate the derivatives
stencil $\{ (\iV-1,\iV,\iV+1,\iV+2), (\jV,\jV+1), (\kV,\kV+1) \}$ with
the cell indices $\ic,\jc,\kc$, and identify $\Usz{1}{xx}(\ic)$,
$\Usz{2}{xx}(\jc)$ and $\Usz{3}{xx}(\kc)$, the core matrices
parametrized with the spatial indices $\ic$, $\jc$, and $\kc$, with
the cores of the second derivative of $\usVTT(\iV,\jV,\kV)$.
Moreover, a comparison with the formulas for the first derivative
shows that:
\begin{align*}
  \usVTTxx(\ic,\jc,\kc) = \frac{ \usVTTx(\ic+1,\jc,\kc) - \usVTTx(\ic-1,\jc,\kc) }{2\hx}.
\end{align*}

Similarly, we have the formulas for the other two second derivatives
\begin{align*}
  \usVTTyy(\ic,\jc,\kc)
  &= \Usz{1}{yy}(\ic)\,\Usz{2}{yy}(\jc)\,\Usz{3}{yy}(\kc)\\[0.5em]
  &:=  
  \frac{1}{2\hy}
  \left( \frac{\Uss{1}(\iV)+\Uss{1}(\iV+1)}{2} \right)
  \left( \frac{\Uss{2}(\jV+2)-\Uss{2}(\jV+1)}{\hy} - \frac{\Uss{2}(\jV)-\Uss{2}(\jV-1)}{\hy} \right)\times\\
  & \hspace{9cm}
  \times\left( \frac{\Uss{3}(\kV)+\Uss{3}(\kV+1)}{2} \right)\\
  &= \frac{ \usVTTy(\ic,\jc+1,\kc) - \usVTTy(\ic,\jc-1,\kc) }{2\hy},
\end{align*}
and
\begin{align*}
  \usVTTzz(\ic,\jc,\kc)
  &= \Usz{1}{zz}(\ic)\,\Usz{2}{zz}(\jc)\,\Usz{3}{zz}(\kc)\\[0.5em]
  &:= 
  \frac{1}{2\hz}
  \left( \frac{\Uss{1}(\iV)+\Uss{1}(\iV+1)}{2} \right)
  \left( \frac{\Uss{2}(\jV)+\Uss{2}(\jV+1)}{2} \right)\times\\
  & \hspace{4cm}
  \times\left( \frac{\Uss{3}(\kV+2)-\Uss{3}(\kV+1)}{\hz} - \frac{\Uss{3}(\kV)-\Uss{3}(\kV-1)}{\hz} \right)\\
  &= \frac{ \usVTTz(\ic,\jc,\kc+1) - \usVTTz(\ic,\jc,\kc-1) }{2\hz}.
\end{align*}
As for the first derivatives, a crucial point in this construction is
that these operations do not introduce any additional error in the TT
representation of the second derivatives other than those of the
finite difference approximation.
Although we will not use them in this work, for completeness' sake, we
provide a possible discretization of the \emph{mixed} second
derivatives of the vertex grid function $\usV$:
\begin{align*}
  \usVTTxy(\ic,\jc,\kc) &=
  \frac{1}{2}\left( 
  \frac{ \usVTTx(\ic,\jc+1,\kc) - \usVTTx(\ic,\jc-1,\kc) }{2\hy}\right.
  \\ & \hspace{5cm}+\left.
  \frac{ \usVTTy(\ic+1,\jc,\kc) - \usVTTy(\ic-1,\jc,\kc) }{2\hx}
  \right),\\[0.5em]
  \usVTTyz(\ic,\jc,\kc) &=
  \frac{1}{2}\left( 
  \frac{ \usVTTy(\ic,\jc,\kc+1) - \usVTTy(\ic,\jc,\kc-1) }{2\hz}\right.
  \\ & \hspace{5cm}+\left.
  \frac{ \usVTTz(\ic,\jc+1,\kc) - \usVTTz(\ic,\jc-1,\kc) }{2\hy} 
  \right),\\[0.5em]
  \usVTTzx(\ic,\jc,\kc) &=
  \frac{1}{2}\left( 
  \frac{ \usVTTz(\ic+1,\jc,\kc) - \usVTTz(\ic-1,\jc,\kc) }{2\hx}\right.
  \\ & \hspace{5cm}+\left.
  \frac{ \usVTTx(\ic,\jc,\kc+1) - \usVTTx(\ic,\jc,\kc-1) }{2\hz}
  \right).
\end{align*}
A straightforward calculation shows that the first formula is
symmetric by swapping $x$ and $y$, i.e., thus formally giving
$\usVTTxy(\ic,\jc,\kc)=\usVTTyx(\ic,\jc,\kc)$;
the second one by swapping $y$ and $z$, thus formally giving
$\usVTTyz(\ic,\jc,\kc)=\usVTTzy(\ic,\jc,\kc)$,
and the third one by swapping $z$ and $x$,  thus formally giving
$\usVTTzx(\ic,\jc,\kc)=\usVTTxz(\ic,\jc,\kc)$.
The efficient implementation must be performed at the core level, thus
leading to formulas similar to the ones for $\usVTTxx$, $\usVTTy$, and
$\usVTTzz$.
However, the average will change the internal ranks and a rounding
step could be necessary.

\begin{remark}
  The alternative formulas that define the discrete second derivatives
  in TT format as the difference between the discrete first
  derivatives in TT format have an obvious correspondence in the
  full-grid formulation, as we can write, for example, that
  \begin{align*}
    \usVxx(\ic,\jc,\kc) = \frac{ \usVx(\ic+1,\jc,\kc) - \usVx(\ic-1,\jc,\kc) }{2\hx}.
  \end{align*}
  Similar difference formulas hold for all the other discrete second
  derivatives.
  Moreover, they provide an alternative way to implement the
  calculation of the second derivatives, which is useful in the next extension of this approach to the case of meshes obtained from
  non-equispaced univariate partitions in $\dirX$, $\dirY$ and $\dirZ$
  and regular equispaced meshes partitioning remapped domains.
\end{remark}
 
\subsection{Discrete Laplacian of vertex grid functions}
\label{subsec:FVM-regular:laplacian}
Summing the discrete second derivatives of the vertex grid function
$\usV$ naturally yields the discrete Laplacian 
\begin{align*}
  \DeltaC\usV(\ic,\jc,\kc) =
  \usVxx(\ic,\jc,\kc) +
  \usVyy(\ic,\jc,\kc) +
  \usVzz(\ic,\jc,\kc). 
\end{align*}
Hereafter, we will refer to the cell grid function $\DeltaC\usV$ as
the \emph{discrete cell Laplacian} of the vertex grid function $\usV$.

\medskip
Since we seek for an approximation of $\us$ at the mesh vertices, we
interpolate $\DeltaC\usV$ at the grid vertices $\V(\iV,\jV,\kV)$ to
obtain the vertex grid function $\DeltaV\usV$, i.e., the \emph{discrete
vertex laplacian} of the grid function $\usV$.
Formally, we introduce a vertex interpolation operator,
denoted by $\IntpV$, such that
\begin{align*}
   \DeltaV\usV = \IntpV\big(\DeltaC\usV\big),
\end{align*}
or, locally, that
\begin{align*}
  \DeltaV\usV(\iV,\jV,\kV) = \IntpV_{\iV,\jV,\kV}\big(\DeltaC\usV\big).
\end{align*}
On a regular mesh where all first neighbor cell-centers around a given
vertex are equidistant from that vertex, and the interpolation is nothing but
the arithmetic average:
\begin{align*}
  \IntpV_{\iV,\jV,\kV}\big(\DeltaC\usV\big)(\iV,\jV,\kV) =
  \frac{1}{8}\sum_{m,n,p=0}^{1}(\DeltaC\usV)(\ic+m,\jc+n,\kc+p).
\end{align*}

\medskip
The same definitions hold for $\usVTT$, the tensor-train
representation of $\usV$, so that the \emph{discrete cell Laplacian in
TT format} reads as
\begin{align*}
  \DeltaCTT\usVTT(\ic,\jc,\kc) =
  \RNDG\big(
  \usVTTxx(\ic,\jc,\kc) +
  \usVTTyy(\ic,\jc,\kc) +
  \usVTTzz(\ic,\jc,\kc)
  \big), 
\end{align*}
where we recall that $\RNDG(\,\cdot\,)$ is the rounding operator that we
discussed at the end of Section~\ref{sec:TTformat}.
The \emph{discrete vertex Laplacian in TT format} reads as
\begin{align*}
   \DeltaVTT\usVTT = \IntpVTT\big(\DeltaCTT\usVTT\big),
\end{align*}
or, locally, as
\begin{align*}
  \DeltaVTT\usVTT(\iV,\jV,\kV) = \IntpVTT_{\iV,\jV,\kV}\big(\DeltaCTT\usVTT\big).
\end{align*}
The vertex interpolation is still the arithmetic average and can be
implemented very efficiently, working directly on the cores of tensor
$\DeltaCTT\usVTT(\ic,\jc,\kc)$:
\begin{align*}
  &\IntpV_{\iV,\jV,\kV}\big(\DeltaCTT\usVTT\big)(\iV,\jV,\kV)
  =
  \frac{1}{8}\sum_{m,n,p=0}^{1}(\DeltaC\usVTT)(\ic+m,\jc+n,\kc+p)\\[0.5em]
  &\qquad\quad
  =
  \left(\frac{\Usz{1}{xx}(\ic)+\Usz{1}{xx}(\ic+1)}{2}\right)
  \left(\frac{\Usz{2}{yy}(\jc)+\Usz{2}{yy}(\jc+1)}{2}\right)
  \left(\frac{\Usz{3}{zz}(\kc)+\Usz{3}{zz}(\kc+1)}{2}\right).
\end{align*}
Such average operation does not change the TT ranks of
$\DeltaCTT\usVTT$, so no additional rounding needs to be performed.

\subsection{Boundary conditions}
\label{subsec:FVM-regular:BCS}
\begin{figure}[htbp]
  \centering
  \begin{tabular}{c c c c}
    \includegraphics[width=.2\linewidth]{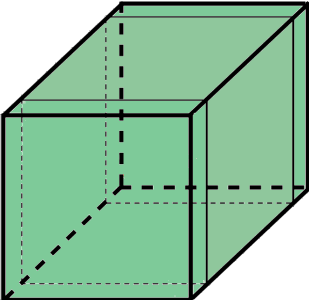} &
    \includegraphics[width=.2\linewidth]{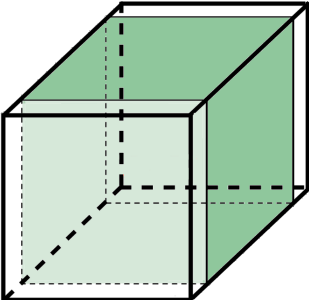} &
    \includegraphics[width=.2\linewidth]{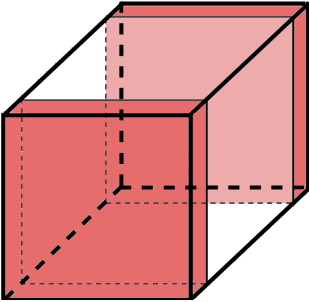} &
    \includegraphics[width=.2\linewidth]{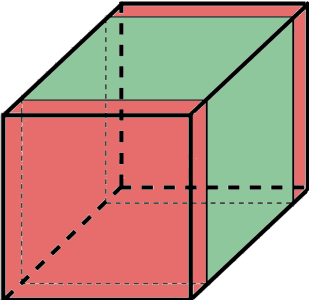} \\
    (a) & (b) & (c) & (d)\\
  \end{tabular}
  \caption{Setting of the boundary conditions along one direction: (a)
    original tensor $\usVTT$; (b) $\usVTT$ with first and last slices
    set to zero; (c) tensor $\tB{}$ with zero inside and exact values on
    the first and last slices; (d) final tensor $\usVTT+\tB{}$.}
  \label{fig:boundary}
\end{figure}

We manage the boundary conditions by using boundary layers of ghost
cells.
In particular, we assume that there are $\Nbnd$ external frames of
cells, which requires $\Nbnd+1$ layers of boundary vertices,
surrounding the hyper-rectangular domain ($\Nbnd=2$ in this work).
These frames correspond to the set of vertex indices:
\begin{itemize}
\item faces orthogonal to $\dirX$: $0\leq\iV\leq\Nbnd$ and
  $\NVx-1-\Nbnd\leq\iV\leq\NVx-1$ for all values of $\jV$ and $\kV$;
\item faces orthogonal to $\dirY$: $0\leq\jV\leq\Nbnd$ and
  $\NVy-1-\Nbnd\leq\jV\leq\NVy-1$ for all values of $\iV$ and $\kV$;
\item faces orthogonal to $\dirZ$: $0\leq\kV\leq\Nbnd$ and
  $\NVz-1-\Nbnd\leq\kV\leq\NVz-1$ for all values of $\kV$ and $\iV$,
\end{itemize}
and to the set of cell indices:
\begin{itemize}
\item faces orthogonal to $\dirX$: $0\leq\ic<\Nbnd$ and
  $\NCx-1-\Nbnd\leq\ic<\NCx-1$ for all values of $\jc$ and $\kc$;
\item faces orthogonal to $\dirY$: $0\leq\jc<\Nbnd$ and
  $\NCy-1-\Nbnd\leq\jc<\NCy-1$ for all values of $\ic$ and $\kc$;
\item faces orthogonal to $\dirZ$: $0\leq\kc<\Nbnd$ and
  $\NCz-1-\Nbnd\leq\kc<\NCx-1$ for all values of $\ic$ and $\jc$.
\end{itemize}
Setting Dirichlet boundary conditions in the full-grid implementation is
straightforward, basically consisting in a loop over all ghost
vertices and cell centers, where the values of the exact solution is
computed and then assigned to $\usV$.
Setting Dirichlet boundary conditions on the TT representations of $\usV$ is
different and a bit trickier since the vertex values of $\usVTT$ are
not stored in the computer memory and are accessible only in an
indirect way.
To enforce such boundary conditions, we modify the TT cores directly, following the approach of
\cite{Manzini-Truong-Vuchkov-Alexandrov:2023}.
The process involves three steps, illustrated in Figure~\ref{fig:boundary}:
\begin{itemize}
\item setting to zero the first and last $\Nbnd$ slices of the
  $\usVTT$ core, corresponding to the boundary vertices;
\item building a new tensor $\tB{}$ whose
  values at the internal indices are zero and whose boundary slices in the cores of $\usVTT$ corresponding to the boundary vertices match a
  decomposition of the boundary data;
\item replacing $\usVTT$ with the sum $\usVTT+\tB{}$.
\end{itemize}  
The boundary slices can be built directly if the boundary condition
functions are 
already in a separated variable form, or can be obtained by performing
a dyadic decomposition, for example, by using the SVD.
If the number of dimensions is bigger than three, the boundary slices
can be obtained by the \TTSVD{} algorithm \cite{oseledets2011} or the
cross-interpolation algorithm \cite{oseledets2010a}.


\section{Extensions to variable-sized grids and remapped domains}
\label{sec:FVM-remapped}

In this section, we extend the finite difference/finite volume discretization framework of the previous section to handle more general computational
domains.
Specifically, we consider two important generalizations: non-uniform
grid spacing with variable-sized cells and domains that can be
remapped onto a cubic domain through orthogonal coordinate transformations.
These extensions allow our numerical methods to be applied to a broader class of problems while maintaining the advantageous
properties of the tensor-train format.

\subsection{Full-grid and tensor-train discretizations on variable-sized grids}
\label{subsec:FVM-remapped:variable_sized}

\begin{figure}[htbp]
\centering
\begin{tabular}{c c}
    \includegraphics[width=.25\linewidth]{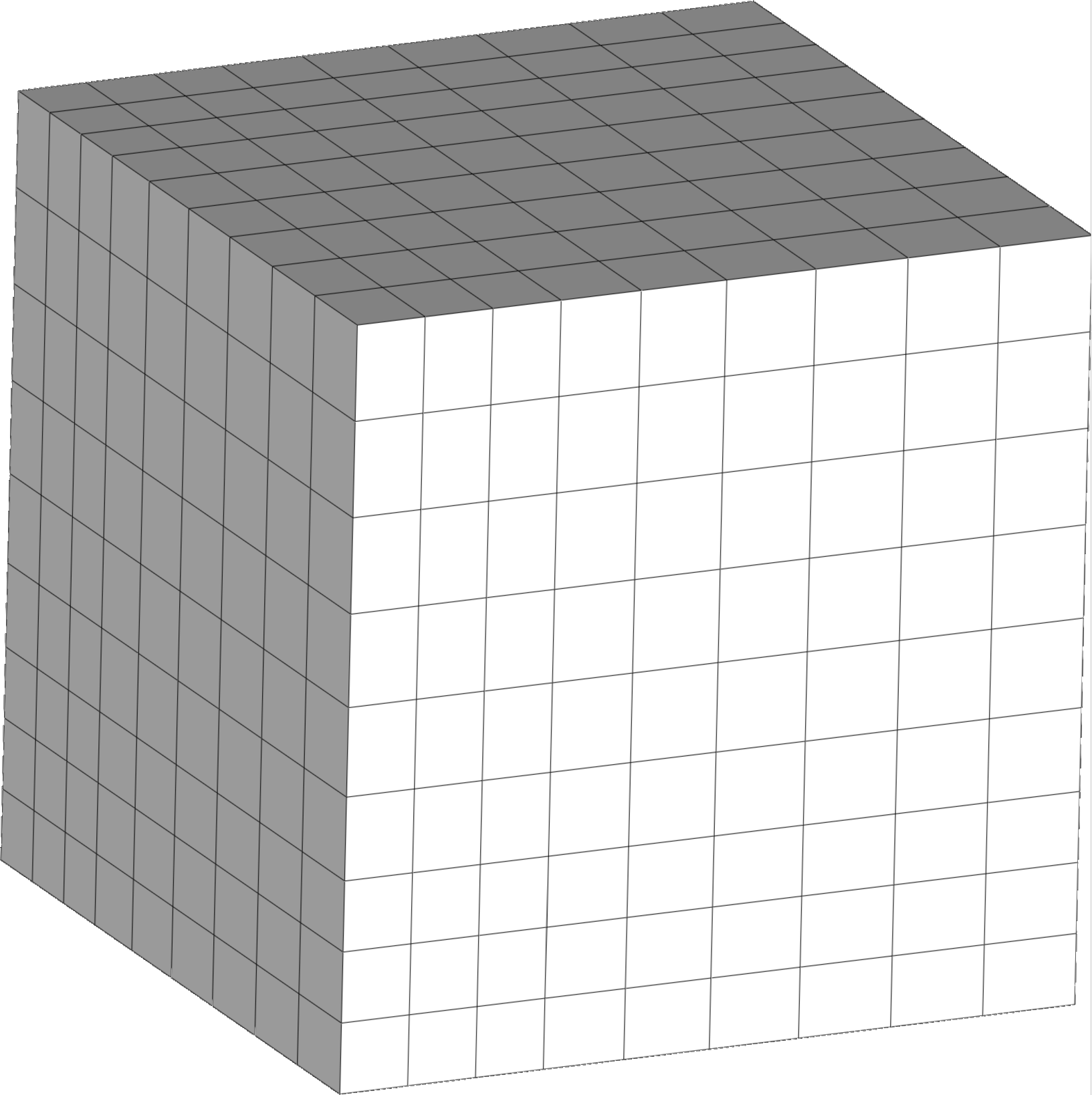} & \qquad
    \includegraphics[width=.25\linewidth]{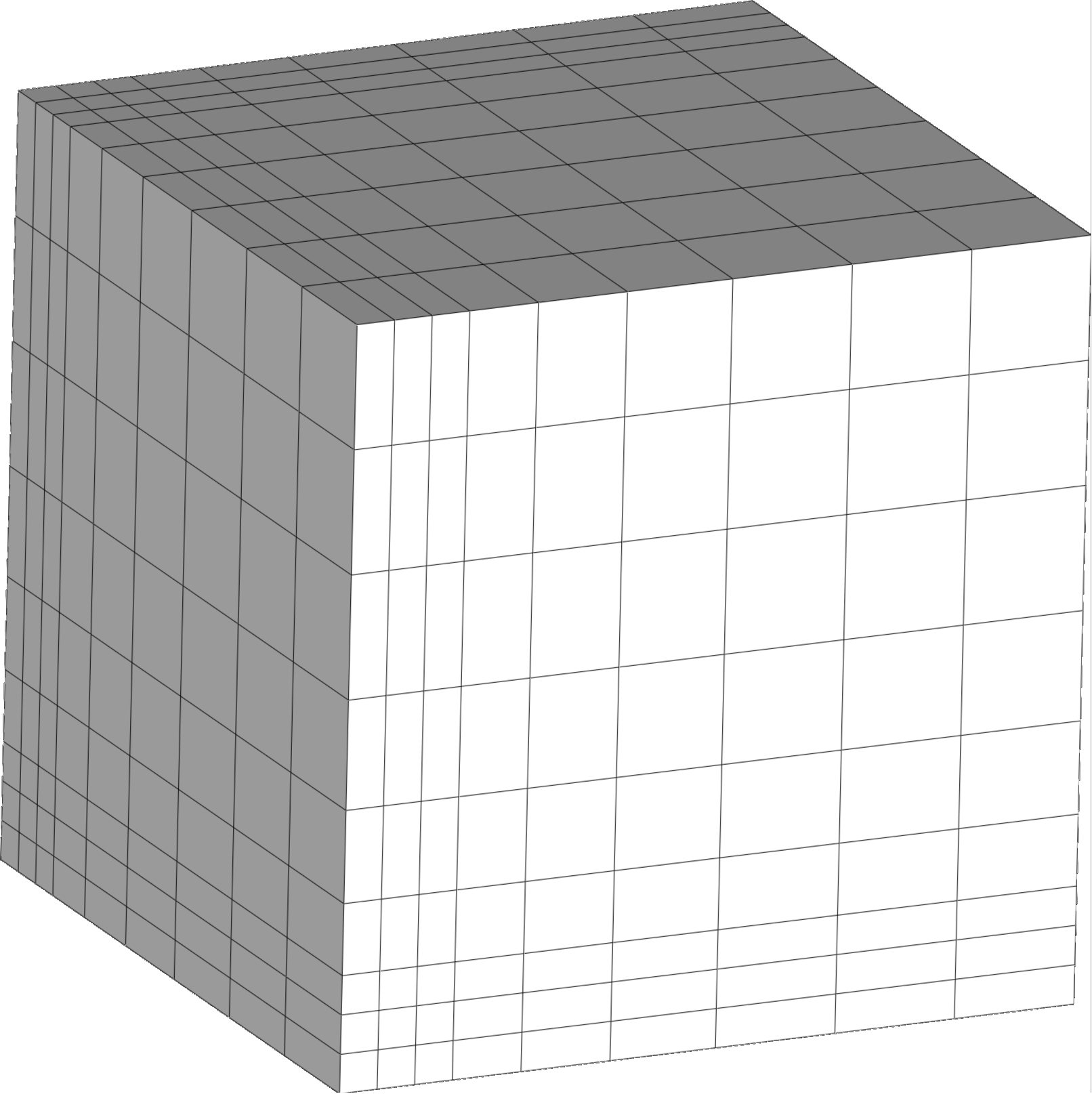}\\
    (a) & (b)\\
\end{tabular}
\caption{Variable-sized grids with geometrically varying step sizes
  $\hxss{\ic}=\rss{\scalebox{0.65}{\dirX}}\hxss{\ic-1}$,
  $\hyss{\jc}=\rss{\scalebox{0.65}{\dirY}}\hyss{\jc-1}$,
  $\hzss{\kc}=\rss{\scalebox{0.65}{\dirZ}}\hzss{\kc-1}$
  and uniform step size ratios:
  $(a)$ $\rss{\scalebox{0.65}{\dirX}}=\rss{\scalebox{0.65}{\dirY}}=\rss{\scalebox{0.65}{\dirZ}}=1.125$;
  $(b)$ $\rss{\scalebox{0.65}{\dirX}}=\rss{\scalebox{0.65}{\dirY}}=\rss{\scalebox{0.65}{\dirZ}}=1.5$}.
\label{fig:grids:variable}
\end{figure}

Assume that the three univariate partitions
$\{\xs(\iV),\,\iV=0,1,\ldots,\NVx-1\}$,
$\{\ys(\jV),\,\jV=0,1,\ldots,\NVy-1\}$, and
$\{\zs(\kV),\,\kV=0,1,\ldots,\NVz-1\}$ along, respectively, $\dirX$,
$\dirY$ and $\dirZ$,
have non-constant step sizes
$\hxss{\ic}=\xs(\iV+1)-\xs(\iV)$,
$\hyss{\jc}=\ys(\jV+1)-\ys(\jV)$, and
$\hzss{\kc}=\zs(\kV+1)-\zs(\kV)$.
The resulting mesh will be an orthogonal grid with a non-uniform grid
spacing, as the one presented in Figure~\ref{fig:grids:variable}.
This type of computational mesh is useful in applications requiring
more accuracy near a domain boundary, for example, in the case the problem solution presents a
boundary layer with a sharp gradient at that boundary.

To obtain the formulas defining the discrete gradient components of a
vertex grid function $\usV$, i.e., $(\usVx,\usVy,\usVz)$ and their TT
representations $(\usVTTx,\usVTTy,\usVTTz)$, we need to use the
variable step sizes $\hxss{\iV}$, $\hyss{\jV}$, and $\hzss{\kV}$ consistently.
With a simple modification, the discrete first derivatives along $\dirX$ are now given
by:
\begin{align*}
  \usVx(\ic,\jc,\kc)   &= \frac{1}{4\hxss{\ic}}\sum_{n=0}^{1}\sum_{p=0}^{1}\Big(\usV(\iV+1,\jV+n,\kV+p) - \usV(\iV,\jV+n,\kV+p)\Big),\\[0.5em]
  \usVTTx(\ic,\jc,\kc) &=
  \Bigg(\frac{\Uss{1}(\iV+1)-\Uss{1}(\iV  )}{\hxss{\ic}}\Bigg)\,
  \Bigg(\frac{\Uss{2}(\jV  )+\Uss{2}(\jV+1)}{2}  \Bigg)\,
  \Bigg(\frac{\Uss{3}(\kV  )+\Uss{3}(\kV+1)}{2}  \Bigg),
\end{align*}
and the formulas for $\usVy$, $\usVz$, $\usVTTy$, and $\usVTTz$ are
obtained accordingly.
The second derivatives are defined similarly.
First, we add and subtract $\usVx(\ic,\jc,\kc)$ to the difference
between the discrete first derivative of $\usV$ at cells
$(\ic+1,\jc,\kc)$ and $(\ic-1,\jc,\kc)$:
\begin{multline*}
  \usVx(\ic+1,\jc,\kc) - \usVx(\ic-1,\jc,\kc) =
  \big[ \usVx(\ic+1,\jc,\kc) - \usVx(\ic,\jc,\kc) \big] \\
  + \big[ \usVx(\ic,\jc,\kc) - \usVx(\ic-1,\jc,\kc) \big].
\end{multline*}
Then, we subdivide the first difference by $\hxss{\ic+1}$ and the
second by $\hxss{\ic}$ and we take the arithmetic average; we find
that
\begin{align*}
  \usVxx(\ic,\jc,\kc) =
  \frac12\frac{ \usVx(\ic+1,\jc,\kc) - \usVx(\ic,  \jc,\kc) }{\hxss{\ic+1}} +
  \frac12\frac{ \usVx(\ic,  \jc,\kc) - \usVx(\ic-1,\jc,\kc) }{\hxss{\ic}}.
\end{align*}
By repeating the same construction, we find the second derivatives
along the other directions:
\begin{align*}
  \usVyy(\ic,\jc,\kc) &=
  \frac12\frac{ \usVy(\ic,\jc+1,\kc) - \usVy(\ic,\jc,  \kc) }{\hyss{\jc+1}} +
  \frac12\frac{ \usVy(\ic,\jc,  \kc) - \usVy(\ic,\jc-1,\kc) }{\hyss{\jc}},\\[0.5em]
  \usVzz(\ic,\jc,\kc) &=
  \frac12\frac{ \usVz(\ic,\jc,\kc+1) - \usVz(\ic,\jc,\kc  ) }{\hzss{\kc+1}} +
  \frac12\frac{ \usVz(\ic,\jc,\kc  ) - \usVz(\ic,\jc,\kc-1) }{\hzss{\kc}},
\end{align*}
and in the TT format:
\begin{align*}
  \usVTTxx(\ic,\jc,\kc) &=
  \frac12\frac{ \usVTTx(\ic+1,\jc,\kc) - \usVTTx(\ic,  \jc,\kc) }{\hxss{\ic+1}} +
  \frac12\frac{ \usVTTx(\ic,  \jc,\kc) - \usVTTx(\ic-1,\jc,\kc) }{\hxss{\ic}},\\[0.5em]
  \usVTTyy(\ic,\jc,\kc) &=
  \frac12\frac{ \usVTTy(\ic,\jc+1,\kc) - \usVTTy(\ic,\jc,  \kc) }{\hyss{\jc+1}} +
  \frac12\frac{ \usVTTy(\ic,\jc,  \kc) - \usVTTy(\ic,\jc-1,\kc) }{\hyss{\jc}},\\[0.5em]
  \usVTTzz(\ic,\jc,\kc) &=
  \frac12\frac{ \usVTTz(\ic,\jc,\kc+1) - \usVTTz(\ic,\jc,\kc  ) }{\hzss{\kc+1}} +
  \frac12\frac{ \usVTTz(\ic,\jc,\kc  ) - \usVTTz(\ic,\jc,\kc-1) }{\hzss{\kc}}.
\end{align*}

\medskip
\noindent
The vertex interpolation algorithm has to be modified accordingly.
For simplicity, we consider the 1D case where we interpolate
cell-centered quantities like $\{\usC(\ic)\}$ that are defined at the
cell centers $\xsC(\ic)$ to the vertex-centered quantity $\usV(\iV)$
at vertex positions $\xsV(\iV)$.
The interpolated value $\usV(\iV)$ is the weighted average on the
stencil $\{\ic-1,\ic\}$:
\begin{align}
  \usV(\iV)
  = \frac{ \usC(\ic-1)\Delta\xss{\ic-1} + \usC(\ic)\Delta\xss{\ic} }{\Delta\xss{\ic-1}+\Delta\xss{\ic}},
  \label{eq:vertex:interpolation:1D}
\end{align}
where $\Delta\xss{\ic-1}=\xsV(\iV)-\xsC(\ic-1)$ is the distance
between vertex $\iV$ and the left cell center, and
$\Delta\xss{\ic}=\xsC(\ic)-\xsV(\iV)$ is the distance between vertex
$\iV$ and the right cell center.
The vertex interpolation in 3D is done by applying the 1D vertex
interpolation dimension by dimension.
For a vertex grid tensor in TT format, we apply the 1D
interpolation independently to each TT core, the resulting
multi-dimensional interpolation being expressed as a tensor product of
the 1D interpolations.

\medskip
\noindent
The construction of the Laplace operator and the setting of the boundary conditions then continues as in
Sections~\ref{subsec:FVM-regular:laplacian}-\ref{subsec:FVM-regular:BCS}.

\subsection{Full-grid and tensor-train discretizations on remapped domains}
\label{subsec:FVM-remapped:remapped_domains}

\begin{figure}[htbp]
\centering
\begin{tabular}{c c}
    \includegraphics[width=.25\linewidth]{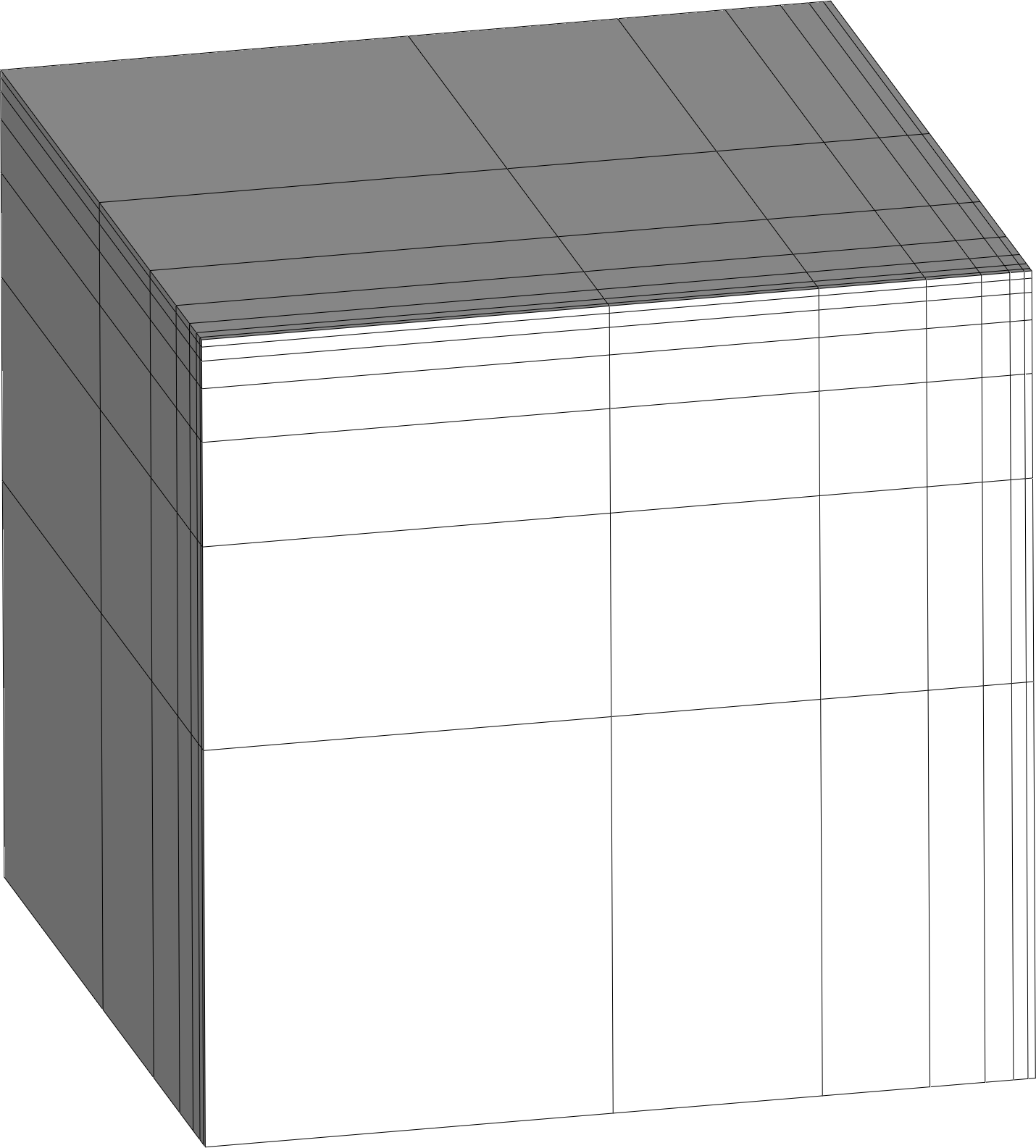} & \qquad
    \includegraphics[width=.25\linewidth]{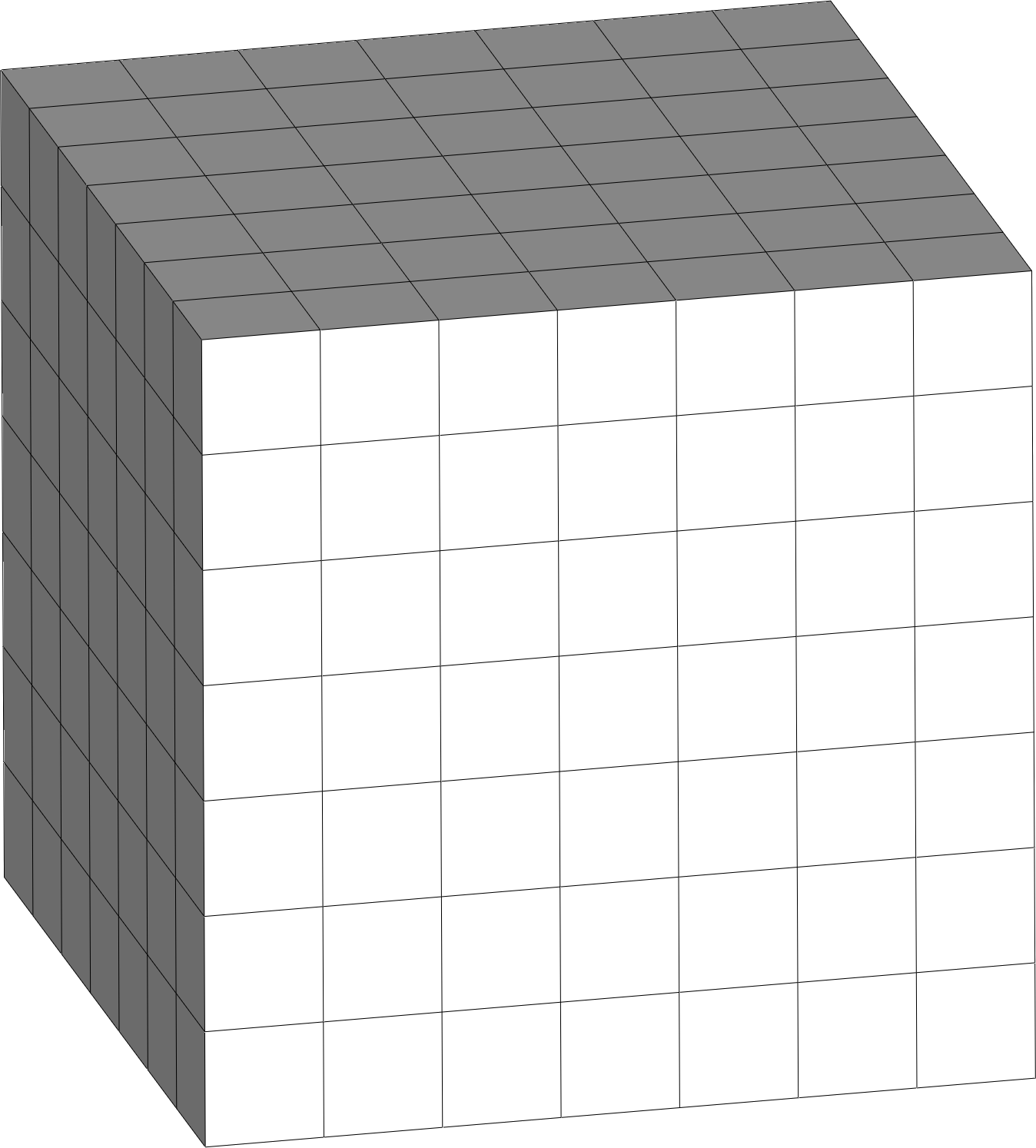}\\
    (a) & (b)\\
\end{tabular}
\caption{Hyper-rectangular 3D domain $\Omega$ (a) and its remapping onto a cube $\widetilde{\Omega}$ (b).}
\label{fig:grids:remapped}
\end{figure}

In the previous sections, we discussed the finite difference
discretization of first and second derivatives in a full-grid and tensor-train format using orthogonal, Cartesian meshes that are
obtained by a tensor product of univariate partitions along
the coordinate axis $\dirX$, $\dirY$, and $\dirZ$.
In this section, we generalize the grid partitions to the case of
hyper-rectangular 3D domains remapped onto a cube, see Figure~\ref{fig:grids:remapped}.
In particular, we assume that the remapping functions are
also univariate so that the resulting multiplicative coefficients do
not affect the rank of the tensor-train finite difference formulas.

Consider the Laplace operator of the function
$\ust:\widetilde{\Omega}\to\REAL$ with respect to the coordinate
framework $(\xst,\yst,\zst)$ defined as:
\begin{align*}
  \widetilde{\Delta}\ust = 
  \frac{\partial^2\ust}{\partial\xst^2} +
  \frac{\partial^2\ust}{\partial\yst^2} +
  \frac{\partial^2\ust}{\partial\zst^2}.
\end{align*}
We establish a mapping between the domains $\Omega$ and
$\widetilde{\Omega}$ through three independent, nonsingular, smooth
coordinate transformations
\begin{align}
  \begin{cases}
    \xs = \xs(\xst),\\
    \ys = \ys(\yst),\\
    \zs = \zs(\zst),\\
  \end{cases}
  \textrm{and their inverse transformations} \quad
  \begin{cases}
    \xst = \xst(\xs), \\
    \yst = \yst(\ys), \\
    \zst = \zst(\zs).
  \end{cases}
  \label{eq:remapping:functions}
\end{align}

\medskip
Consider the function $\us:\Omega\to\REAL$ that corresponds to the
function $\ust$ defined on $\widetilde{\Omega}$ through the coordinate
map introduced above so that
$\us(\xs,\ys,\zs)=\ust(\xst(\xs),\yst(\ys),\zst(\zs))$.
Its Laplace operator is given by:
\begin{align}
  \Delta\us &= 
  \left(\frac{\partial\xst(\xs)}{\partial\xs}\right)^{-1}\frac{\partial}{\partial\xs}\left[
    \left(\frac{\partial\xst(\xs)}{\partial\xs}\right)^{-1}\frac{\partial\us}{\partial\xs}\right]
  +
  \left(\frac{\partial\yst(\ys)}{\partial\ys}\right)^{-1}\frac{\partial}{\partial\ys}\left[ 
    \left(\frac{\partial\yst(\ys)}{\partial\ys}\right)^{-1}\frac{\partial\us}{\partial\ys}\right] +\nonumber\\[0.5em]
  &+
  \left(\frac{\partial\zst(\zs)}{\partial\zs}\right)^{-1}\frac{\partial}{\partial\zs}\left[ 
    \left(\frac{\partial\zst(\zs)}{\partial\zs}\right)^{-1}\frac{\partial\us}{\partial\zs}\right].
  \label{eq:remapped:Laplacian}
\end{align}

The construction of discrete derivatives proceeds in the same way as
for variable-sized grids.
Since the heat equation can be discretized on the remapped
domain by using a regular, equispaced, orthogonal Cartesian grid, we assume 
$\hx$, $\hy$ and $\hz$ to be constant.
However, we must take into account the derivatives of the maps
$\xst(\xs)$, $\yst(\ys)$, and $\zst(\zs)$ that appear in the remapped
equations.
We obtain the formulas defining the discrete gradient components of a
vertex grid functions $\usV$, i.e., $(\usVx,\usVy,\usVz)$ and its TT
representation $(\usVTTx,\usVTTy,\usVTTz)$, by consistently
multiplying the step sizes $\hx$, $\hy$ and $\hz$ by
$\xst'(\ic)=\partial\xst(\xsC(\ic))\slash{\partial\xs}$,
$\yst'(\jc)=\partial\yst(\ysC(\jc))\slash{\partial\ys}$, and
$\zst'(\kc)=\partial\zst(\zsC(\kc))\slash{\partial\zs}$.
Therefore, the discrete first derivatives along $\dirX$, $\dirY$ and $\dirZ$
are given by:
\begin{align*}
  \usVx(\ic,\jc,\kc)   &=
  \frac{1}{4\xst'(\xs(\ic))\hx}\sum_{n=0}^{1}\sum_{p=0}^{1}\Big(\usV(\iV+1,\jV+n,\kV+p) - \usV(\iV,\jV+n,\kV+p)\Big)\\[0.5em]
  \usVy(\ic,\jc,\kc)   &=
  \frac{1}{4\yst'(\ys(\jc))\hy}\sum_{m=0}^{1}\sum_{p=0}^{1}\Big(\usV(\iV+m,\jV+1,\kV+p) - \usV(\iV+m,\jV,\kV+p)\Big),\\
  \usVz(\ic,\jc,\kc)   &=
  \frac{1}{4\zst'(\zs(\kc))\hz}\sum_{m=0}^{1}\sum_{n=0}^{1}\Big(\usV(\iV+m,\jV+n,\kV+1) - \usV(\iV+m,\jV+n,\kV)\Big),
\end{align*}
and their TT version are
\begin{align*}
  \usVTTx(\ic,\jc,\kc) &=
  \Bigg(\frac{\Uss{1}(\iV+1)-\Uss{1}(\iV  )}{\xst'(\xsC(\ic))\hx}\Bigg)\,
  \Bigg(\frac{\Uss{2}(\jV  )+\Uss{2}(\jV+1)}{2}  \Bigg)\,
  \Bigg(\frac{\Uss{3}(\kV  )+\Uss{3}(\kV+1)}{2}  \Bigg),\\
  \usVTTy(\ic,\jc,\kc) &=
  \Bigg(\frac{\Uss{1}(\iV  )+\Uss{1}(\iV+1)}{2}  \Bigg)\,
  \Bigg(\frac{\Uss{2}(\jV+1)-\Uss{2}(\jV  )}{\yst'(\ysC(\jc))\hy}\Bigg)\,
  \Bigg(\frac{\Uss{3}(\kV  )+\Uss{3}(\kV+1)}{2}  \Bigg),\\
  \usVTTz(\ic,\jc,\kc) &=
  \Bigg(\frac{\Uss{1}(\iV  )+\Uss{1}(\iV+1)}{2}  \Bigg)\,
  \Bigg(\frac{\Uss{2}(\jV  )+\Uss{2}(\jV+1)}{2}  \Bigg)\,
  \Bigg(\frac{\Uss{3}(\kV+1)-\Uss{3}(\kV  )}{\zst'(\zsC(\kc))\hz}\Bigg).
\end{align*}

The second derivatives are obtained by repeating the same derivation
as for variable-sized grids, still taking into account that now each mesh size is a constant value over the mesh but multiplied by a variable function from the remapping process.
We obtain:
\begin{align*}
  \usVxx(\ic,\jc,\kc) &=
  \frac12\frac{ \usVx(\ic+1,\jc,\kc) - \usVx(\ic,  \jc,\kc) }{\xst'(\xsC(\ic+1))\hx} +
  \frac12\frac{ \usVx(\ic,  \jc,\kc) - \usVx(\ic-1,\jc,\kc) }{\xst'(\xsC(\ic))\hx},\\[0.5em]
  \usVyy(\ic,\jc,\kc) &=
  \frac12\frac{ \usVy(\ic,\jc+1,\kc) - \usVy(\ic,\jc,  \kc) }{\yst'(\ysC(\jc+1))\hy} +
  \frac12\frac{ \usVy(\ic,\jc,  \kc) - \usVy(\ic,\jc-1,\kc) }{\yst'(\ysC(\jc))\hy},\\[0.5em]
  \usVzz(\ic,\jc,\kc) &=
  \frac12\frac{ \usVz(\ic,\jc,\kc+1) - \usVz(\ic,\jc,\kc  ) }{\zst'(\zsC(\kc+1))\hz} +
  \frac12\frac{ \usVz(\ic,\jc,\kc  ) - \usVz(\ic,\jc,\kc-1) }{\zst'(\zsC(\kc))\hz},
\end{align*}
and in the TT format:
\begin{align*}
  \usVTTxx(\ic,\jc,\kc) &=
  \frac12\frac{ \usVTTx(\ic+1,\jc,\kc) - \usVTTx(\ic,  \jc,\kc) }{\xst'(\xsC(\ic+1))\hx} +
  \frac12\frac{ \usVTTx(\ic,  \jc,\kc) - \usVTTx(\ic-1,\jc,\kc) }{\xst'(\xsC(\ic))\hx},\\[0.5em]
  \usVTTyy(\ic,\jc,\kc) &=
  \frac12\frac{ \usVTTy(\ic,\jc+1,\kc) - \usVTTy(\ic,\jc,  \kc) }{\yst'(\ysC(\jc+1))\hy} +
  \frac12\frac{ \usVTTy(\ic,\jc,  \kc) - \usVTTy(\ic,\jc-1,\kc) }{\yst'(\ysC(\jc))\hy},\\[0.5em]
  \usVTTzz(\ic,\jc,\kc) &=
  \frac12\frac{ \usVTTz(\ic,\jc,\kc+1) - \usVTTz(\ic,\jc,\kc  ) }{\zst'(\zsC(\kc+1))\hz} +
  \frac12\frac{ \usVTTz(\ic,\jc,\kc  ) - \usVTTz(\ic,\jc,\kc-1) }{\zst'(\zsC(\kc))\hz}.
\end{align*}

\medskip
\noindent
The vertex interpolation is done by applying the same algorithm of
\eqref{eq:vertex:interpolation:1D} as for the case of variable size
grids, the only difference being in the definition of the
interpolation coefficients:
\begin{itemize}
\item $\Delta\xss{\ic-1}=\big(\xst'(\xsV(\iV))-\xst'(\xsC(\ic-1))\big)\hx$
  is the distance between vertex $\V(\iV)$ and the center of
  $\C(\ic-1)$, the cell on the left;
\item $\Delta\xss{\ic}=\big(\xst'(\xsC(\ic))-\xst'(\xsV(\iV))\big)\hx$
  is the distance between vertex $\V(\iV)$ and the center of
  $\C(\ic)$, the cell on the right.
\end{itemize}

\medskip
\noindent
The construction of the Laplace operator and the setting of the
boundary conditions then continue as in
Sections~\ref{subsec:FVM-regular:laplacian}-\ref{subsec:FVM-regular:BCS}.
\begin{remark}
  The variable-sized grids of
  Section~\ref{subsec:FVM-remapped:variable_sized} can be expressed
  through a similar remapping algorithm when the vertex positions are
  determined by a set of remapping functions as
  in~\eqref{eq:remapping:functions}.
  In such a case, the numerical treatment is the same.
  However, the case considered
  in~Section~\ref{subsec:FVM-remapped:variable_sized} is more general
  as it does not require an explicit knowledge of such remapping
  functions.
  This situation occurs, for example, when the mesh vertices are
  adaptively relocated in the domain according to an a-posteriori
  error indicator.
\end{remark}



\section{Time integration}
\label{sec:time_integration}

We partition the time interval $[0,\Ts]$ using a time step with
constant size $\Delta\ts$, and denote the intermediate time instants
as $\ts^n=n\Delta\ts$ for $n = 0,1,\ldots,\NT$, where
$\NT=\Ts/\Delta\ts$ is the total number of time steps.
Let $\usV^n$ and $\usV^{\TT,n}$ be the time discretizations
$\usV(\cdot;\ts^n)$ and $\usVTT(\cdot;\ts^n)$, the full-grid and
tensor-train vertex grid functions that approximate $\us(\cdot,\ts)$
at $\ts=\ts^n$.
We discretize the time derivatives of $\usV(\cdot;t)$ and
$\usVTT(\cdot;t)$ using the first-order accurate in time finite
difference formula, so that
$\big(\partial\usV\slash{\partial\ts}\big)^n\approx\big(\usV^{n+1}-\usV^{n}\big)\slash{\Delta\ts}$
and
$\big(\partial\usVTT\slash{\partial\ts}\big)^n\approx\big(\usV^{\TT,n+1}-\usV^{\TT,n}\big)\slash{\Delta\ts}$.

\medskip
To update $\usV^{n}$ to $\usV^{n+1}$ and $\usV^{\TT,n}$ to
$\usV^{\TT,n+1}$, we consider three different integration schemes: the
explicit Euler scheme, the implicit Euler scheme, and the semi-implicit
Crank-Nicolson scheme.
The explicit Euler scheme, while straightforward to implement, is only
first-order accurate in time, and its simplicity comes at the cost of
a potential instability that may require smaller time steps.
It is worth noting that the explicit and implicit Euler schemes are
also building blocks of higher-order algorithms, e.g., the Runge-Kutta
scheme, with improved accuracy and bigger storage and computational
costs.
In our formulation, we adopt a similar strategy and split the
Crank-Nicolson scheme in two steps with an initial explicit step
followed by an implicit one, both advancing the numerical solution of
$\Delta\ts/2$.

\medskip
Let $\fsV^{n+\theta}$, $\theta\in\{0,1/2,1\}$, denote the source term
evaluated at the time $\ts^{n}+\theta\Delta\ts$ at the mesh vertices.
The TT reformulation requires a TT decomposition of this term, which
can be performed using a tensor-train decomposition algorithm such as
the \TTSVD{} or the \emph{cross-interpolation algorithm}, see, e.g.,
\cite{oseledets2011,oseledets2010a}, and the comment at the end of
Section~\ref{sec:TTformat}.
We let $\fsVTT(\cdot;\ts)$ denote the tensor-train representation of
the vertex grid function $\fsV(\cdot;\ts)$ that is obtained by
sampling the forcing function $\fs(\cdot,\ts)$ at time $\ts$ at the
mesh vertices.
For $\theta=1/2$, as in the Crank-Nicolson scheme, we use the
approximation $\fs^{n+1/2}\approx(\fs^{n}+\fs^{\ns+1})/2$, which is
still second-order accurate and does not affect the global accuracy of
our method.

\medskip
Hereafter, we will refer to the vertex grid functions
$\Delta\ts\big(\DeltaV\usV+\fsV\big)$ and the corresponding tensor-train 
representation $\Delta\ts\big(\DeltaVTT\usVTT+\fsVTT\big)$ as
the \emph{vertex residual terms} or, simply, \emph{residual term}.
Finally, we note that the construction of both $\DeltaV\usV$ and
$\DeltaVTT\usVTT$ already takes into account the Dirichlet boundary
conditions, which are explicitly set in the ghost vertex and cell
frames at the same time instant at which we consider the fields $\usV$
and $\usVTT$.
Therefore, we only perform the update of $\usV$ and $\usVTT$ at the
internal vertices.

\paragraph{Explicit Euler scheme} The explicit Euler scheme
evaluates the residual term at the time $\ts^n$, so that
\begin{align}
  \usV^{n+1} = \usV^{n} + \Delta\ts\big( \DeltaV\usV^{n} + \fsV^{n} \big).
  \label{eq:explicit}
\end{align}
In the TT reformulation, we first perform a similar update:
\begin{align}
  \usVb^{\TT,n+1} = \usV^{\TT,n} + \Delta\ts\big( \DeltaVTT\usV^{\TT,n} + \fsV^{\TT,n} \big),
  \label{eq:explicit:TT:update}
\end{align}
again for all index triple $(\iV,\jV,\kV)$ corresponding to the
\emph{internal} vertices, and then we apply the rounding algorithm \begin{align}
  \usV^{\TT,n+1} = \RNDG\big(\usVb^{\TT,n+1}\big).
  \label{eq:explicit:TT:rounding}
\end{align}
Algorithm~\ref{algo:explicit:TT} in~\ref{sec:appx:time_integration} provides
implementation details on the update from $\usV^{\TT,n}$ to
$\usV^{\TT,n+1}$.

\paragraph{Implicit Euler scheme} The implicit Euler method
evaluates the residual term at the time $\ts^{n+1}$, so that
\begin{align}
  \usV^{n+1} - \Delta\ts\DeltaV\usV^{n+1} = \usV^{n} + \Delta\ts\fsV^{n+1}.
  \label{eq:implicit}
\end{align}
In the TT reformulation, we perform a similar update:
\begin{align}
  \usV^{\TT,n+1} - \Delta\ts\DeltaV\usV^{\TT,n+1} = \usV^{\TT,n} + \Delta\ts\fsV^{\TT,n+1}.
  \label{eq:implicit:TT}
\end{align}
Both formulations \eqref{eq:implicit} and~\eqref{eq:implicit:TT}
require the resolution of a linear system that we iteratively solve by
applying a \emph{matrix-free Preconditioned Conjugate Gradient (PCG)
method}.
Krylov methods need to perform the \emph{matrix-vector multiplication}
of the coefficient matrix times the current conjugate direction at
every iteration.
In the matrix-free setting, we perform this by evaluating the
left-hand side of~\eqref{eq:implicit} and~\eqref{eq:implicit:TT} on
the conjugate direction fields.
The details of the implementation and the design of a preconditioner
that is suitable to the TT formulation are discussed in the next
subsection.
Algorithm~\ref{algo:implicit:TT} in~\ref{sec:appx:time_integration}
provides implementation details on the update from $\usV^{\TT,n}$ to
$\usV^{\TT,n+1}$, including the rounding procedure that is performed
inside the PCG to control the rank growth.

\paragraph{Semi-implicit Crank-Nicolson scheme} The
semi-implicit Crank-Nicolson scheme evaluates the residual term at the
half-time step $\ts^{n+1/2}=\ts+\Delta\ts/2$.
We approximate this evaluation by taking the average of the residual
at the initial and final times, i.e., $\ts^{n}$ and $\ts^{n+1}$ so
that the time marching schemes read as:
\begin{align*}
  \usV^{n+1} - \frac{\Delta\ts}{2}\big( \DeltaV\usV^{n+1} + \fsV^{n+1} \big) =
  \usV^{n}  + \frac{\Delta\ts}{2}\big( \DeltaV\usV^{n} + \fsV^{n} \big),
\end{align*}
and
\begin{align*}
  \usV^{\TT,n+1} - \frac{\Delta\ts}{2}\big( \DeltaVTT\usV^{\TT,n+1} + \fsV^{\TT,n+1} \big) =
  \usV^{\TT,n}  + \frac{\Delta\ts}{2}\big( \DeltaVTT\usV^{\TT,n} + \fsV^{\TT,n} \big).
\end{align*}
This approach allows us to split the time step into two consecutive
half steps.
The first step uses the explicit Euler scheme to advance the solution
from $\usV^{n}$ $\big(\text{or }\usV^{\TT,n}\big)$ to an intermediate solution
$\usV^{n+1/2}$ $\big(\text{or }\usV^{\TT,n+1/2}\big)$ at time
$\ts^{n+1/2}=\ts^{n}+\Delta\ts/2$.
The second step uses the implicit Euler scheme to advance the
intermediate solution from $\usV^{n+1/2}$ $\big(\text{or }\usV^{\TT,n+1/2}\big)$ to
the final time step $\usV^{n+1}$ $\big(\text{or }\usV^{\TT,n+1}\big)$.
We can thus reformulate the Crank-Nicolson scheme as follows: 
\begin{itemize}
\item[$(1)$] \textit{Explicit step:} use formula~\eqref{eq:explicit}:
  \vspace{-1.5\lineskip}
  \begin{align*}
    \usVb^{\ns+1/2}    &= \usV^{n}  + \frac{\Delta\ts}{2}\big( \DeltaV\usV^{n} + \fsV^{n} \big),
    \intertext{\vspace{-1\lineskip} or formula \eqref{eq:explicit:TT:update}:}
    \usVb^{\TT,\ns+1/2} &= \usV^{\TT,n}  + \frac{\Delta\ts}{2}\big( \DeltaVTT\usV^{\TT,n} + \fsV^{\TT,n} \big),
  \end{align*}
  with $\Delta\ts/2$ instead of $\Delta\ts$;
  
\item[$(2)$] \textit{Implicit step:} use formula~\eqref{eq:implicit}:
  \vspace{-1.5\lineskip}
  \begin{align*}
    \usV^{n+1}    - \frac{\Delta\ts}{2}\big( \DeltaV\usV^{n+1} + \fsV^{n+1} \big)         &= \usVb^{\ns+1/2},\\
    \intertext{\vspace{-0.75\lineskip} or formula \eqref{eq:implicit:TT}:}
    \usV^{\TT,n+1} - \frac{\Delta\ts}{2}\big( \DeltaVTT\usV^{\TT,n+1} + \fsV^{\TT,n+1} \big) &= \usVb^{\TT,\ns+1/2},\\[-0.5em]
  \end{align*}
  with $\Delta\ts/2$ instead of $\Delta\ts$.
\end{itemize}
It is worth mentioning that the boundary conditions are set at the
initial time $\ts^{n}$ in the first explicit step and at time
$\ts^{n+1/2}$ in the second and last implicit step.
Algorithm~\ref{algo:Crank-Nicolson:TT}
in~\ref{sec:appx:time_integration} provides implementation details on
the update from $\usV^{\TT,n}$ to $\usV^{\TT,n+1}$.

\subsection{Matrix-free preconditioned conjugate gradient method}
\label{subsec:time_integration:PCG}

The PCG routine is detailed in Algorithm~\ref{algo:pcg}.
This algorithm is an implementation of the Conjugate Gradient
  method, adapted to operate on vectors in the tensor-train format as
  described below and with the introduction of a rounding step for the
  rank growth control.
  We propose a possible preconditioning strategy for this algorithm,
  which applies independently to each TT core.
  Our thorough experimentation indicates that this approach performs
  effectively, and we present it as the practical method we employed
  for solving the linear systems arising from implicit and
  semi-implicit time-marching schemes.
  A theoretical investigation of this formulation within the framework
  of Krylov subspace methods is beyond the scope of this paper and
  will be addressed in future research.

We discuss here only the implicit/semi-implicit step for the
tensor-train, as the full-grid case is similar.
We initialize the algorithm with the field $\usV^{\TT,n+\tau}$, where
$\tau=1$ for the implicit Euler scheme, in which we use $\Delta\ts$,
and $\tau=1/2$ for the semi-implicit Crank-Nicolson scheme, where we
use $\Delta\ts/2$ instead.
The algorithm needs two external procedures, \MATVECPROD{} and
\PRECON, and applies the rounding algorithm after every update of the
solution $\usVTT$, the auxiliary residual $\rsVTT$, and the conjugate
direction $\psVTT$.
The procedure \MATVECPROD{} applies the discrete heat operator to an
input TT field $\psVTT$, i.e.,
$\psVTT\to(\psVTT-\tau\Delta\ts\DeltaVTT\psVTT)$, where $\DeltaVTT$ is
the discrete vertex Laplace operator defined in
Section~\ref{sec:FVM-regular}.
  
\begin{algorithm}[hbt!]
    \caption{Preconditioned Conjugate Gradient (PCG)}
    \label{algo:pcg}
    \textbf{Input:}
    \begin{itemize}[nosep]
    \item $\usV^{\TT,n+\tau}$: vertex grid function at time $t^{n+\tau}$, with $\tau\in\left\{1/2,1\right\}$
    \item $\bsVTT$: right-hand side function on $\Omega$
    \item $\text{tol}$: tolerance
    \item $\text{maxiter}$: maximum number of iterations
    \end{itemize}
    \textbf{Output:}
    \begin{itemize}[nosep]
    \item $\usV^{\TT,n+1}$: vertex grid function at time $t^{n+1}$
    \end{itemize}
  \begin{algorithmic}[1]
    \Procedure{PCG}{$\usV^{\TT,n+\tau}, \bsVTT, \text{tol}, \text{maxiter}$}
    \State $\usVTT \gets \usV^{\TT,n+\tau}$                                          \Comment{Initial guess for the solution}
    \State $\rsVTT \gets \RNDG(\bsVTT - \MATVECPROD(\usVTT))$                     \Comment{Initial residual}
    \State $\|\rsV^{\TT,0}\| \gets \|\rsVTT\|$                                 \Comment{Compute \& store the initial residual norm}
    \State $\ztVTT \gets \PRECON(\rsVTT)$                                          \Comment{Apply preconditioner}
    \State $\psVTT \gets \ztVTT$                                                   \Comment{Assign Initial search direction}
    \State $\rho \gets \rsVTT\cdot\ztVTT$
    \For{$\ell \gets 1$ \textbf{to} $\text{maxiter}$}                              \Comment{Krylov loop}
    \State $(\As\psVTT) \gets \MATVECPROD(\psVTT)$                                 \Comment{Matrix-vector product}
    \State $\alpha \gets \rho / (\psVTT\cdot(\As\psVTT))$                          \Comment{Update $\alpha$}
    \State $\usVTT \gets \RNDG(\usVTT + \alpha \psVTT)$                            \Comment{Update solution}
    \State $\rsVTT \gets \RNDG(\rsVTT - \alpha (\As\psVTT))$                       \Comment{Update residual}
    \If {$\|\rsVTT\| < \text{tol}\,\|\rsV^{\TT,0}\|$}                             \Comment{Check convergence}
    \State $\usV^{\TT,n+1} \gets \usVTT$
    \State \Return $\usV^{\TT,n+1}$                                                  \Comment{Converged!}
    \EndIf
    \State $\ztVTT \gets \PRECON(\rsVTT)$                                          \Comment{Apply preconditioner}
    \State $\rho_{\text{new}} \gets \rsVTT\cdot\ztVTT$                                \Comment{Update $\rho$}
    \State $\psVTT \gets \RNDG(\ztVTT + (\rho_{\text{new}} / \rho) \psVTT)$           \Comment{Update conjugate direction}
    \State $\rho \gets \rho_{\text{new}}$
    \EndFor
    \State $\usV^{\TT,n+1} \gets \usVTT$  \Comment{Maximum iterations reached}
    \EndProcedure
  \end{algorithmic}
\end{algorithm}

For this operator, we construct the preconditioner $\P$ based on the
following splitting of the Laplacian operator:
\begin{align*}
  \Ps =
  \big(\alpha\Is - \tau\Delta\ts\Pss{x}\big)
  \big(\alpha\Is - \tau\Delta\ts\Pss{y}\big)
  \big(\alpha\Is - \tau\Delta\ts\Pss{z}\big),
\end{align*}
see~\cite{Benzi:2002,Ng:2004}.
Here, $\Pss{\ell}=\TRIDIAG(1,2,1)\slash{\hh_{\ell}^2}$
with
$\ell\in\{x,y,z\}$ represents the tridiagonal matrices corresponding to the standard 3-point stencil finite difference discretization of the
univariate Laplacian operator.
This construction effectively approximates the inverse of the heat
operator by a product of one-dimensional operators.
The procedure \PRECON{} applies this preconditioner to the residual
$\rsVTT$, returning the preconditioned vector field
$\zeta^{\TT}_{\V}$.
Applying the preconditioner $\Ps$ to $\rsVTT$ involves solving the
linear system:
\begin{align*}
  \ztVTT
  = \Ps^{-1}\rsVTT
  = \big(\Is - \tau\Delta\ts\Pss{z}\big)^{-1}\big(\Is - \tau\Delta\ts\Pss{y}\big)^{-1}\big(\Is - \tau\Delta\ts\Pss{x}\big)^{-1}\rsVTT.
\end{align*}
Consider the TT representations of $\zeta^{\TT}_{\V}$ and $\rsVTT$
(using the symbols $\tZss{\ell}$ and $\tRss{\ell}$ to denote the
respective cores):
\begin{align*}
  \ztVTT(\iV,\jV,\kV)
  &= \sum_{\alpha_1=1}^{\rss{1}}\sum_{\alpha_2=1}^{\rss{2}}
  \tZss{1}(\iV,\alpha_1)\tZss{2}(\alpha_1,\jV,\alpha_2)\tZss{3}(\alpha_2,\kV),\\[0.5em]
  \rsVTT(\iV,\jV,\kV)
  &= \sum_{\alpha_1=1}^{\rss{1}}\sum_{\alpha_2=1}^{\rss{2}}
  \tRss{1}(\iV,\alpha_1)\tRss{2}(\alpha_1,\jV,\alpha_2)\tRss{3}(\alpha_2,\kV),
\end{align*}
where $\ztVTT$ and $\rsVTT$ have the same ranks $\rss{1}$ and
$\rss{2}$.
The application of $\Ps$ reduces to solving a series of tridiagonal
systems along each spatial dimension defined by the matrices
$\Pss{x}$, $\Pss{y}$, and $\Pss{z}$ and the corresponding spatial
fibers of $\rsVTT$.
The cores of $\ztVTT$ are then given by:
\begin{align*}
  &
  \forall\alpha_1=1,2,\ldots,\rss{1},\quad
  \forall\alpha_2=1,2,\ldots,\rss{2}
  \\[0.5em] &\qquad
  \begin{array}{rlcll}
    &\tZss{1}(         :,\alpha_1) &=& \big(\Is - \tau\Delta\ts\Pss{x}\big)^{-1}\tRss{1}(         :,\alpha_1),\\[0.5em]
    &\tZss{2}(\alpha_1,:,\alpha_2) &=& \big(\Is - \tau\Delta\ts\Pss{y}\big)^{-1}\tRss{2}(\alpha_1,:,\alpha_2),\\[0.5em]
    &\tZss{3}(\alpha_2,:         ) &=& \big(\Is - \tau\Delta\ts\Pss{z}\big)^{-1}\tRss{3}(\alpha_2,:         ).
  \end{array}
\end{align*}
These tridiagonal systems can be efficiently solved using the Thomas
algorithm~\cite{Varga:2001}.



\section{Numerical Experiments}
\label{sec:numerical}

This section presents a numerical investigation of the tensor-train
method developed in the previous sections applied to
problem~\eqref{eq:problem:strong}-\eqref{eq:problem:strong2}.
To this end, we first verify the consistency of our Laplace operator
approximation to confirm that the discretization error scales as
$\calO(\hh^2)$ as expected when refining the mesh, where
$\hx=\hy=\hz:=\hs$.
After this validation, we investigate the rank stability of the
numerical approximate solution when
problem~\eqref{eq:problem:strong}-\eqref{eq:problem:strong2} admits
rank-1 solution.
Finally, we evaluate the method's effectiveness using a manufactured
solution, varying the initial approximation's rank to assess both the
computational efficiency and accuracy of the proposed method.

To impose that $\hx=\hy=\hz$, we set $\NCx=\NCy=\NCz:=\NC$ in all test
cases.
We consider a sequence of grid meshes defined on the unit cube
$\Omega=[0,1]^3$, starting from a $20\times20\times20$ grid ($\NC=20$
and $\hs=1/\NC=0.05$) with two additional boundary layers ($\Nbnd=2$).
We compare the FG and the TT implementations in three different scenarios:
\begin{itemize}
\item \textit{regular grids:} tensor product grids with equispaced
  distributed vertices in each direction;
\item \textit{variable grids:} tensor product grids with variable mesh
  size in each direction, such as in Figure~\ref{fig:grids:variable}.
  In particular, $\hs$ degrades by a factor of 1.125 in each
  direction, so that $\hxss{0}=\hs$ and
  $\hxss{\ic}=1.125\ \hxss{\ic-1}$ for each $\ic$, the same holding
  for $\hyss{\jc}$ and $\hzss{\kc}$;
\item \textit{remapped domains:} tensor product of equispaced
  distributed vertices on a remapped domain, as the one presented in
  Figure~\ref{fig:grids:remapped}.
  The transformation gives the remapping:
  \begin{align*}
    \left(\xst, \yst, \zst \right) \to \left(-2e^{-2\xs}, -2e^{-2\ys}, -2e^{-2\zs}\right).
  \end{align*}
\end{itemize}
All errors are evaluated using the Frobenius norm suitably rescaled by
$\hh^{\frac32}$, to have a norm that is equivalent to a mesh-dependent
$\LTWO$-like norm.
In this section, we will label all the results from the full-grid
implementations by ``\emph{FG}'' and all the results from the
tensor-train implementations by ``\emph{TT}''.
Computations have been performed using a 3.60GHz Intel Core i7
processor with 16 cores and 32GB of RAM.

\subsection{Consistency error evaluation}
\label{subsec:numerical:consistency}
Let $\us(\xv)$ be a smooth enough function (at least
$\CS{2}$-regular), $\usV$ the grid function that is obtained by
sampling $\us$ at the grid nodes, and $\usVTT$ its tensor-train
representation.
We define the consistency error of the Laplace operator approximation
at $\us$ as the residual norms
\begin{align*}
  \RESDV(\usV)   = \norm{ \DeltaV\usV     - \restrict{\big(\Delta\us\big)}{\V} }{}
  \qquad\textrm{and}\qquad
  \RESDV(\usVTT) = \norm{ \DeltaVTT\usVTT - \restrict{\big(\Delta\us\big)}{\V} }{},  
\end{align*}
where
$\restrict{\big(\Delta\us\big)}{\V}=\big\{\Delta\us(\xsV,\ysV,\zsV)\big\}$,
the Laplacian of $\us$ evaluated at the grid nodes $\V$.
The consistency error is expected to scale as $\calO(\hh^2)$ when
$\hh$ is refined.
We report the results for the approximation of the Laplacian of
\begin{align}
  \uspace(\xs,\ys,\zs)=\sin(2\pi\xs)\,\sin(2\pi\ys)\,\sin(2\pi\zs),\quad\xs,\ys,\zs\in[0,1],
  \label{eq:space}
\end{align}
in Table~\ref{tab:consistency}.
We measure the residual norms $\RESD_{\V,FG}(\usV)$,
$\RESD_{\V,TT}(\usVTT)$ relative to the full-grid and the tensor-train
implementations, and their convergence rates, denoted as
$\rate_{FG},\ \rate_{TT}$.
We also report computational times ($\ctime_{FG},\ \ctime_{TT}$) and
storage ($\cstrg_{FG},\ \cstrg_{TT}$), with ratios
$\ctime_r=\ctime_{FG}/\ctime_{TT}$ and
$\cstrg_r=\cstrg_{FG}/\cstrg_{TT}$.

\begin{table}[htbp]
\caption{Consistency error tests: residual $\RESDV$, convergence rate
  $\rate$, computational time $\ctime$, and storage $\cstrg$ in FG and
  TT implementation, with ratios $\ctime_r=\ctime_{FG}/\ctime_{TT}$
  and $\cstrg_r=\cstrg_{FG}/\cstrg_{TT}$.  }
\label{tab:consistency}
\centering
\begin{subtable}[h]{.9\textwidth}
  \caption{Regular grids.}
  \begin{tabular}{c}
    \includegraphics[width=\textwidth]{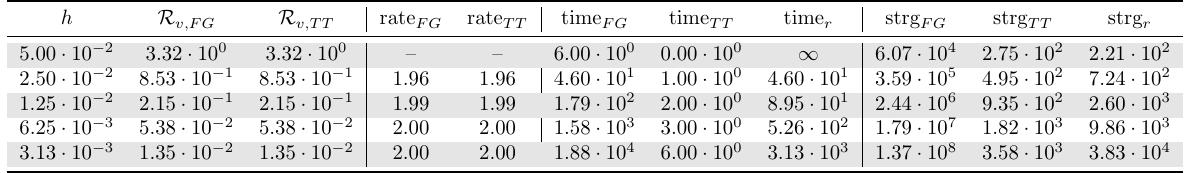}
  \end{tabular}
\end{subtable}
\begin{subtable}[h]{.9\textwidth}
  \caption{Variable sized grids.}
  \begin{tabular}{c}
    \includegraphics[width=\textwidth]{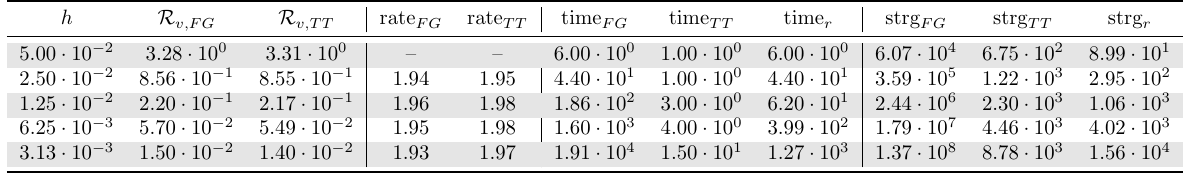}
  \end{tabular}
\end{subtable}
\begin{subtable}[h]{.9\textwidth}
  \caption{Remapped domains.}
  \begin{tabular}{c}
    \includegraphics[width=\textwidth]{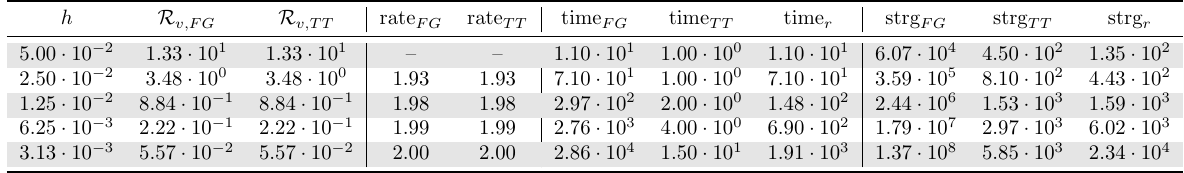}
  \end{tabular}
\end{subtable}
\end{table}

Both methods achieve almost the same consistency error levels and
exhibit the expected second-order convergence rate in all the test
cases, thus confirming the consistency of the TT Laplace operator.
We also note that the TT solver outperforms the FG solver in storage
requirements and computational time by multiple orders of magnitude,
demonstrating the extraordinary efficiency of the TT approach.

\subsection{Rank-1 solution approximations}
\label{subsec:numerical:rank1}
As a second experiment, we consider the rank-1 and time-dependent solution
\begin{align}
  \utime(x,y,z,t)=\sin(2\pi x)\ \sin(2\pi y)\ \sin(2\pi z)\ \cos(2\pi t) \quad x,y,z,t\in[0,1],
  \label{eq:time}
\end{align}
comparing the explicit Euler, implicit Euler, and Crank-Nicolson
schemes with an initial time step equal to $\Delta\ts=10^{-4}$.
In this case, at each refinement level, in addition to halving the
space size $\hs$, we divide $\Delta\ts$ by a coefficient 4 in the
explicit Euler method (for stability reasons) and by a coefficient 2
in the implicit Euler and semi-implicit Crank-Nicolson methods.
The number of time steps $\NT$ is set to $1/\Delta\ts$, ensuring that
the simulation always reaches the final time $\Ts=1$.
For computational reasons, in the last two time refinements of the FG
implementation, we set $\NT=100$ and multiply the computational time by
$1/(\Delta\ts\ \NT)$ to estimate the total computational time.
We empirically observed that this approximation leads to an average
$3\%$ discrepancy compared with the actual computational time.
The rounding threshold was set to $\toll=10^{-4}$; this constant
affects the $\RNDG$ procedure present in the setting of the boundary
conditions, in PCG and in the update of the solution at each time
step.

\medskip
The results for all tests are reported in
Tables~\ref{tab:regular}-\ref{tab:variable}-\ref{tab:remapped}.
In this case, and from now on, we measure the relative approximation
errors $\varepsilon_{FG}=\norm{\usV-\us}{}/\norm{\us}{}$ and
$\varepsilon_{TT}=\norm{\usV^{TT}-\us}{}/\norm{\us}{}$.
Since both the explicit and implicit Euler methods are only
first-order in time, and the semi-implicit method is a second-order
accurate method in time, we expect to see a convergence rate
proportional to $\calO(\hh^2+\Delta\ts)$ in the first case and an
error scaling like $\calO(\hh^2+\Delta\ts^2)$ in the second case.
In the last two rows of each table, we report the estimated time
values in italics, while all storage values are exact as they do not
depend on the number of time steps.

\begin{table}[htbp]
\caption{Rank-1 solution on regular grids: relative error
  $\varepsilon$, convergence rate $\rate$, computational time
  $\ctime$, and storage $\cstrg$ in FG and TT implementation, with
  ratios $\ctime_r=\ctime_{FG}/\ctime_{TT}$ and
  $\cstrg_r=\cstrg_{FG}/\cstrg_{TT}$.}
\label{tab:regular}
\centering
\begin{subtable}[h]{\textwidth}
  \caption{Explicit Euler method.}
  \begin{tabular}{c}
    \includegraphics[width=\textwidth]{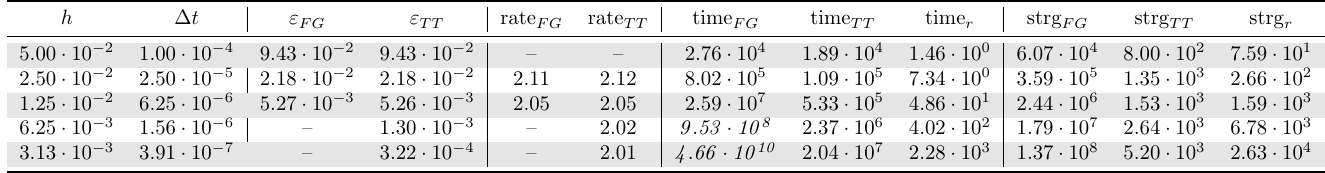}
  \end{tabular}
\end{subtable}
\begin{subtable}[h]{\textwidth}
  \caption{Implicit Euler method.}
  \begin{tabular}{c}
    \includegraphics[width=\textwidth]{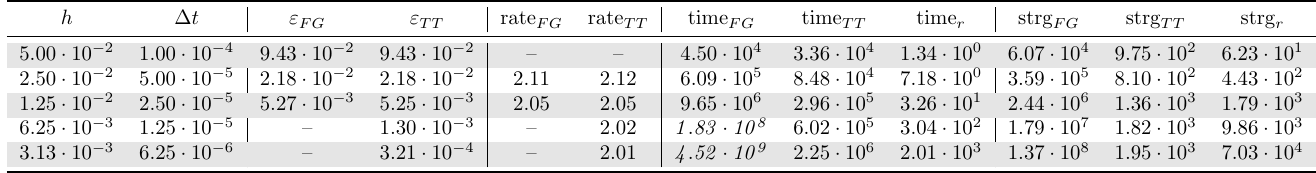}
  \end{tabular}
\end{subtable}
\begin{subtable}[h]{\textwidth}
    \caption{Crank-Nicolson method.}
    \begin{tabular}{c}
      \includegraphics[width=\textwidth]{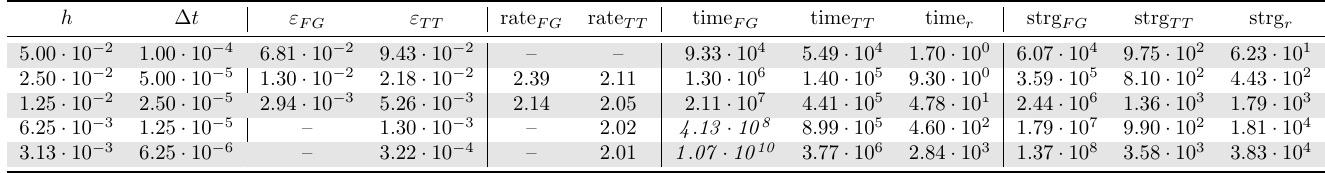}
    \end{tabular}
\end{subtable}
\end{table}

\begin{table}[htbp]
  \caption{Rank-1 solution on variable sized grids: relative error
    $\varepsilon$, convergence rate $\rate$, computational time
    $\ctime$, and storage $\cstrg$ in FG and TT implementation, with
    ratios $\ctime_r=\ctime_{FG}/\ctime_{TT}$ and
    $\cstrg_r=\cstrg_{FG}/\cstrg_{TT}$.}
  \label{tab:variable}
  \centering
  \begin{subtable}[h]{\textwidth}
    \caption{Explicit Euler method.}
    \begin{tabular}{c}
      \includegraphics[width=\textwidth]{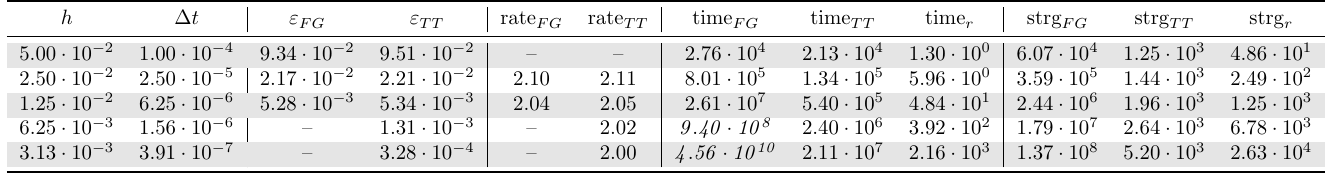}
    \end{tabular}
  \end{subtable}
  \begin{subtable}[h]{\textwidth}
    \caption{Implicit Euler method.}
    \begin{tabular}{c}
      \includegraphics[width=\textwidth]{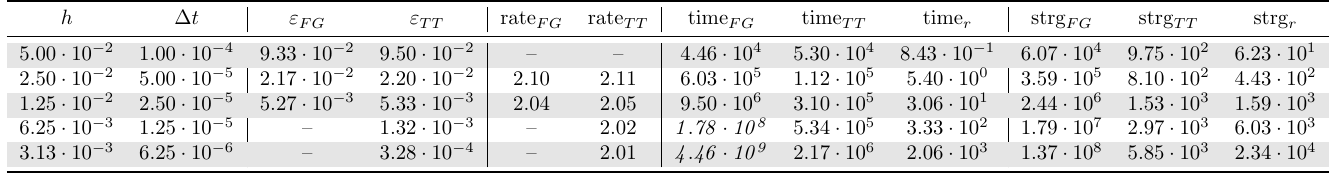}
    \end{tabular}
  \end{subtable}
  \begin{subtable}[h]{\textwidth}
    \caption{Crank-Nicolson method.}
    \begin{tabular}{c}
      \includegraphics[width=\textwidth]{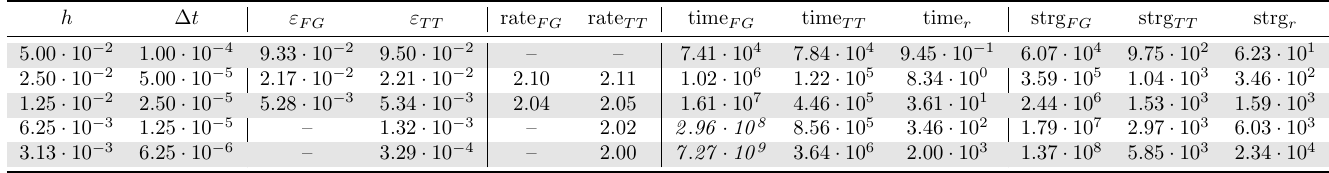}
    \end{tabular}
  \end{subtable}
\end{table}

\begin{table}[htbp]
  \caption{Rank-1 solution on remapped domains: relative error
    $\varepsilon$, convergence rate $\rate$, computational time
    $\ctime$, and storage $\cstrg$ in FG and TT implementation, with
    ratios $\ctime_r=\ctime_{FG}/\ctime_{TT}$ and
    $\cstrg_r=\cstrg_{FG}/\cstrg_{TT}$.}
  \label{tab:remapped}
  \centering
  \begin{subtable}[h]{\textwidth}
    \caption{Explicit Euler method.}
    \begin{tabular}{c}
      \includegraphics[width=\textwidth]{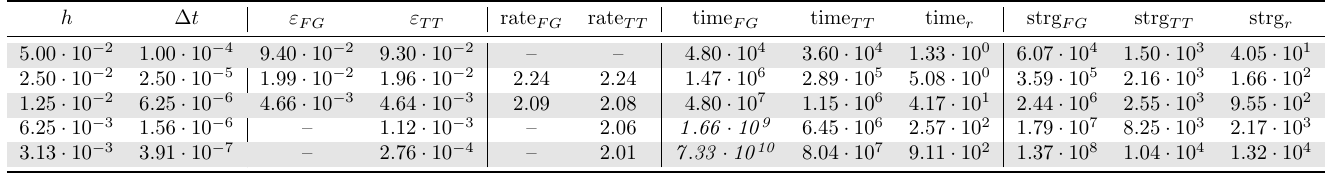}
    \end{tabular}
  \end{subtable}
  \begin{subtable}[h]{\textwidth}
    \caption{Implicit Euler method.}
    \begin{tabular}{c}
      \includegraphics[width=\textwidth]{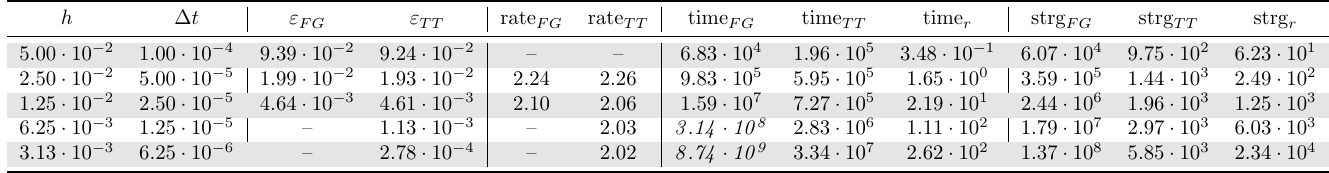}
    \end{tabular}
  \end{subtable}
  \begin{subtable}[h]{\textwidth}
    \caption{Crank-Nicolson method.}
    \begin{tabular}{c}
      \includegraphics[width=\textwidth]{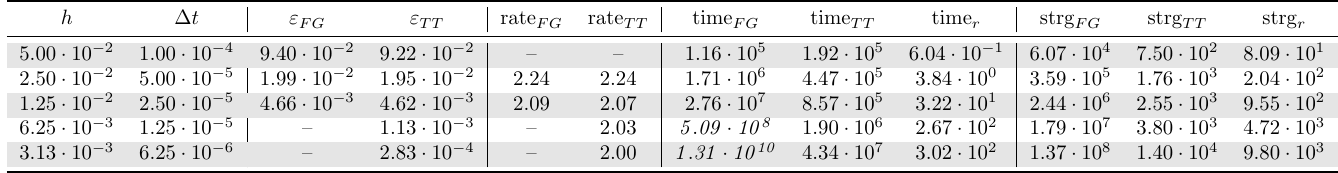}
    \end{tabular}
  \end{subtable}
\end{table}

The tensor-train reformulation demonstrates remarkable efficiency,
while maintaining the same accuracy as the FG approach.
Indeed, both methods exhibit optimal second-order convergence rates,
as expected from the scheme construction, in all test cases, i.e.,
regular grids, variable-size grids along the Cartesian dimensions, and
regular grids on remapped domains, and the three time-marching options
of explicit, implicit, and semi-implicit schemes.
However, the TT solver dramatically outperforms the FG solver in
storage requirements and computational time by multiple orders of
magnitude.
In fact, we can see that:
\begin{itemize}
\item \emph{the computational time ratio} (columns
  $\ctime_r=\ctime_{FG}/\ctime_{TT}$) increases from $\calO(10)$ to
  $\calO(10^3)$, with greater benefits at finer resolutions.
\item \emph{the storage ratio} (columns
  $\cstrg_r=\cstrg_{FG}/\cstrg_{TT}$) increases consistently from
  $\calO(10)$ on the coarser grids to $\calO(10^4)$ on the finer
  grids;
\end{itemize}
Such efficiency gains are achieved without compromising accuracy, as
shown by the error levels in the columns labeled by $\varepsilon_{FG}$
and $\varepsilon_{TT}$, which are identical throughout all refined
grid calculations.

\medskip
As explained in Section~\ref{sec:TTformat}, the maximum TT rank of
$\usVTT$ plays a crucial role in the simulation, related to its
accuracy and computational cost.
Throughout Sections~\ref{sec:FVM-regular}-\ref{sec:FVM-remapped}, we
pointed out how some operations, e.g. the sums of derivative tensors,
while not introducing additional errors in the TT representation of
the finite difference approximation, may increase the rank thus
requiring apposite rounding steps to keep our tensor-train
representation ``low-rank''.
In Figure~\ref{fig:ranks}, we plot the maximum TT rank reached by
$\usVTT$ throughout the whole simulation, relative to the regular grid
scenario with the explicit and implicit Euler method
(cf. Tables~\ref{tab:regular}(a)-(b)).
Each line corresponds to a different space refinement (the rows of
Table~\ref{tab:regular}), and the $x-$values indicate the advancement
in time.
Notably, the ranks do not increase in time, and finer grids lead to
lower ranks; in fact, in the last two refinements, the rank is
constantly equal to one.
We refer to this phenomenon as the ``rank stability'' of the
algorithm, and it indicates that our implementation can keep a
low-rank representation of a low-rank solution during the whole time
integration.
\begin{table}[htpb]
  \centering
  \begin{tabular}{c c}
    \includegraphics[width=0.45\linewidth]{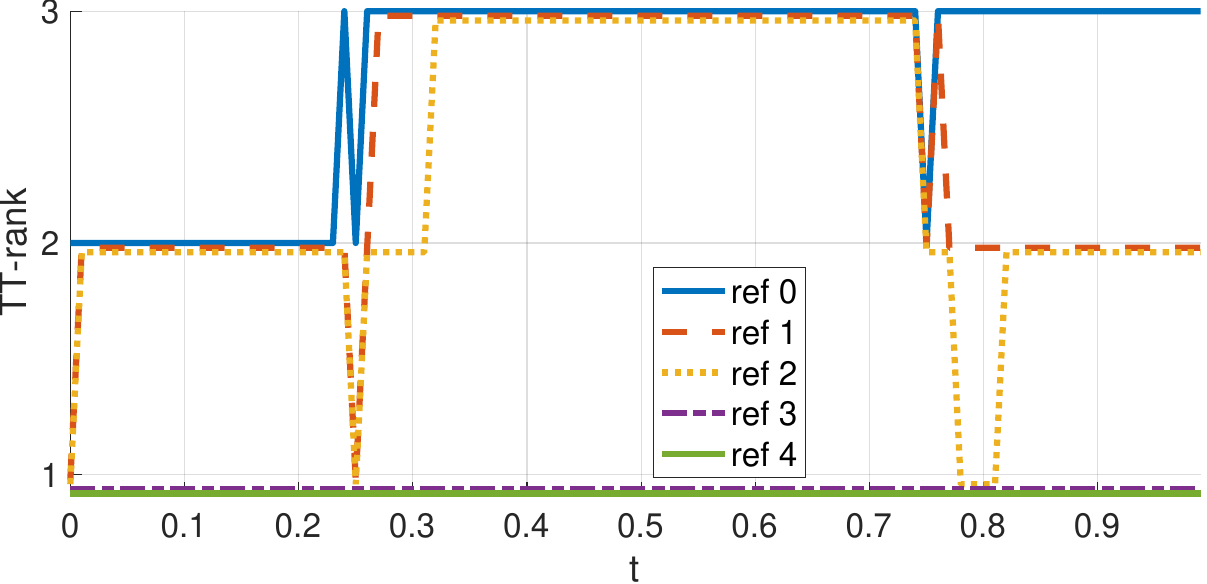} &
    \includegraphics[width=0.45\linewidth]{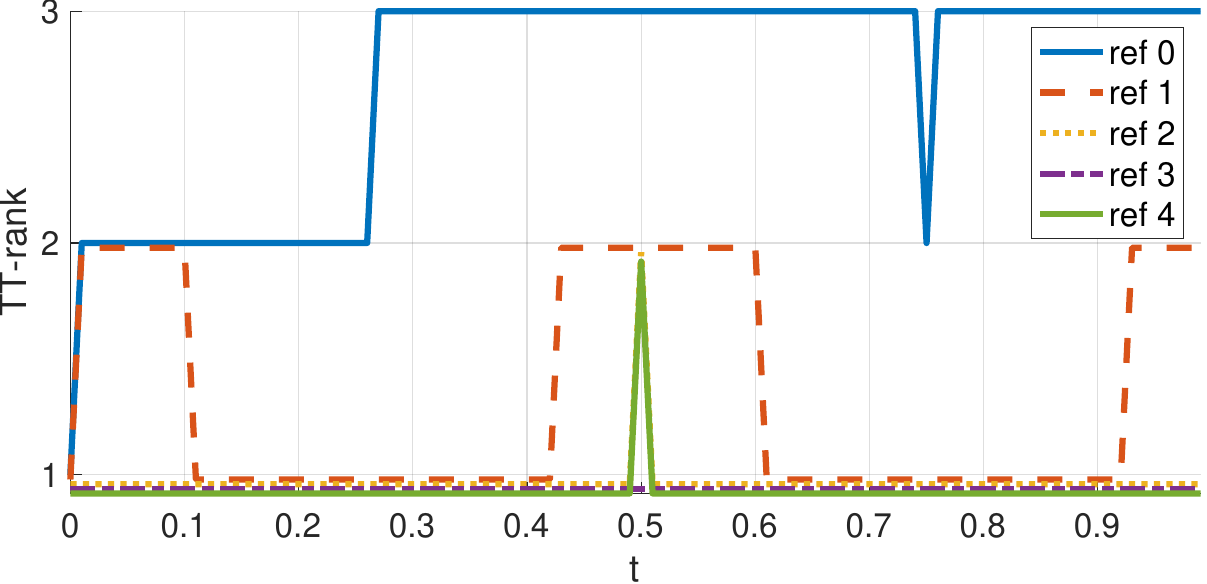} \\
                    {\footnotesize (a) explicit Euler method} & 
                    {\footnotesize (b) implicit Euler method}
  \end{tabular}
  \caption{Maximum rank of $\usVTT$ in each refinement step of the
    Euler method for rank-1 solution on regular grids.  Plot values
    have been slightly shifted vertically to improve readability.}
  \label{fig:ranks}
\end{table}

Last, we investigate the impact of the preconditioner on the condition
number $\kappa(A)=\lambda_{\max}(A)/\lambda_{\min}(A)$ of the global
stiffness ``matrix'' $A=I-\Delta\ts\,\DeltaVTT$ that we use in the
tensor-train implementation of the Euler implicit method (the same
argument holds for the Crank-Nicolson semi-implicit method).
Being $A$ a perturbation of the identity ``matrix'' $I$ for small
$\Delta\ts$, we find that the minimum eigenvalue of $A$ is
$\lambda_{\min}(A)\simeq1$, so that we can (roughly) estimate the
condition numbers $\kappa(A)$ and $\kappa(P^{-1}AP)$ by computing only
$\lambda_{\max}(A)$ and $\lambda_{\max}(P^{-1}AP)$ through the power
method, see~\cite[Chapter~7.3]{Golub-VanLoan:2013}.
Note that the power method only requires the matrix-vector product
operation that we can perform by applying the routine \MATVECPROD{},
cf. Section~\ref{subsec:time_integration:PCG}.
In Table~\ref{tab:condition_number}, we consider the results obtained
for a sequence of refined regular grids for different values of
$\Delta\ts$, e.g., $\Delta\ts\in\{1,0.1,0.01\}$.
We compare the maximum eigenvalue of the ``plain matrix'' $A$ and the
``preconditioned matrix'' $P^{-1}AP$ for decreasing $\hh$ values.
If $\Delta\ts$ is sufficiently small, the stiffness ``matrix'' almost
reduces to the identity ``matrix'' and both $\lambda_{\max}(A)$ and
$\lambda_{\max}(P^{-1}AP)$ remain small, while bigger time steps lead
to a significant growth of the maximum eigenvalue in the
non-preconditioned case.
The PCG effectively addresses this undesirable effect, maintaining
$\lambda_{\max}(P^{-1}AP)\simeq1$.
Importantly, the use of the preconditioner strongly mitigates the
correlation between $\lambda_{\max}(P^{-1}AP)$ and $\hh$, thus
allowing the solver to work easily with much smaller mesh sizes.

\begin{table}[htpb]
  \centering
  \caption{Maximum eigenvalue $\lambda_{\max}$ of the stiffness matrix
    $A$, with and without preconditioner $P$, using different time
    steps $\Delta\ts$.}
  \includegraphics[width=.8\linewidth]{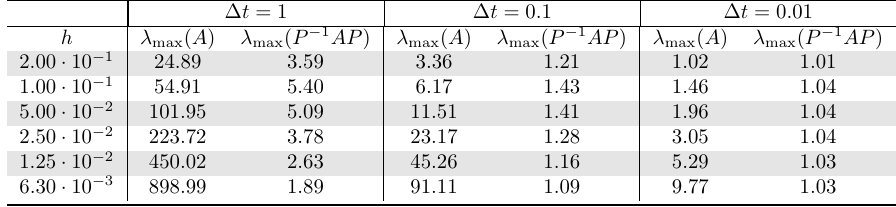}
  \label{tab:condition_number}
\end{table}

\subsection{Rank-stable approximations of a general  solution}
\label{subsec:numerical:high-rank}
As a further experiment, on a sequence of regular grids we consider
the more general solution
\begin{align}
  \urank(x,y,z,t)=e^{-\left(x+y+z-\frac{3}{2}\right)^2} t, \quad x,y,z,t\in[0,1].
\label{eq:rank}
\end{align}
In the case of high-rank solutions, a bad choice of the ranks for the
initial approximations may compromise the accuracy of the entire
simulation.
From an initial calculation, we found that the grid function obtained
by sampling only the space-dependent part of $\urank(x,y,z,t)$ can be
represented in tensor-train format with TT ranks equal to 11 and with
an approximation error of the order of the machine precision.
Then, in a set of preliminary tests (the results of which are not
reported) we discovered that an initial solution approximated with
ranks strictly less than $5$ leads to bad approximations.
We also noticed that having to deal with higher function ranks, the
rounding threshold used in the tests of
Section~\ref{subsec:numerical:rank1} was not small enough to provide
accurate results.
Therefore, we consider the three different initial ranks 5, 11, and
15, with $\toll=10^{-7}$.
In Table~\ref{tab:high-rank}, we report the results for the explicit
Euler method.
As in Section~\ref{subsec:numerical:rank1}, we start with $\hs=0.05$
and $\Delta\ts=10^{-4}$, halving both of them at each refinement
level, with $\NT=1000$ time steps.

\begin{table}[htbp]
  \caption{General solution on regular grids, explicit Euler method:
    relative error $\varepsilon$, convergence rate $\rate$,
    computational time $\ctime$, and storage $\cstrg$ in FG and TT
    implementation, with ratios $\ctime_r=\ctime_{FG}/\ctime_{TT}$ and
    $\cstrg_r=\cstrg_{FG}/\cstrg_{TT}$.}
  \label{tab:high-rank}
  \centering
  \begin{subtable}[h]{\textwidth}
    \caption{Initial rank = 5.}
    \begin{tabular}{c}
      \includegraphics[width=\textwidth]{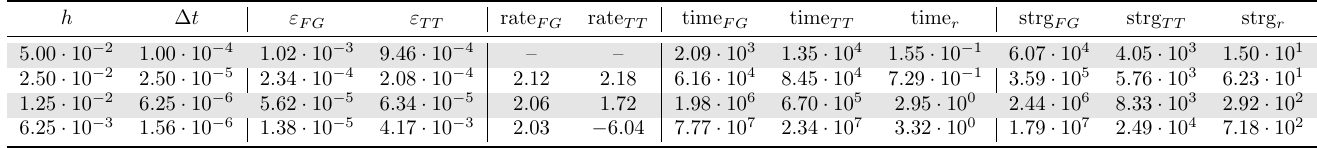}
    \end{tabular}
  \end{subtable}
  \begin{subtable}[h]{\textwidth}
    \caption{Initial rank = 11.}
    \begin{tabular}{c}
      \includegraphics[width=\textwidth]{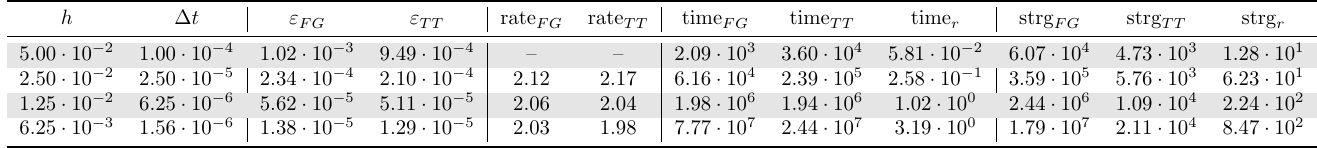}
    \end{tabular}
  \end{subtable}
  \begin{subtable}[h]{\textwidth}
    \caption{Initial rank = 15.}
    \begin{tabular}{c}
      \includegraphics[width=\textwidth]{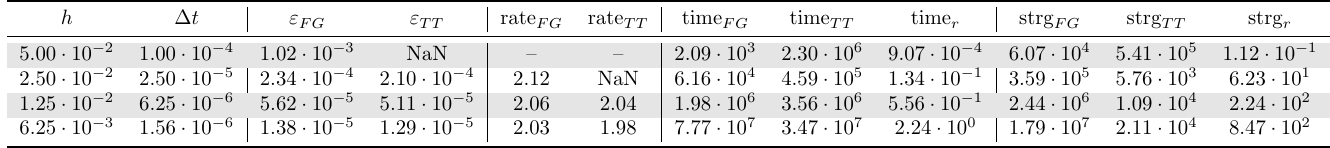}
    \end{tabular}
  \end{subtable}
\end{table}

Using an initial rank of $5$ leads to smaller tensors and, therefore,
faster performance.
However, after some refinements, such a small rank is not enough to
preserve the accuracy required by the scheme, and the convergence
deteriorates.
With an initial rank of $11$ instead, the error decreases as expected,
and the method converges with the correct rate.
After a few space refinements, the computational time ratio becomes
favorable to the TT solver, while the storage required by the TT
implementation is always much smaller compared to the FG.
Last, using a higher initial rank implies more expensive computations
and does not improve the solution accuracy.
In the first line of Table~\ref{tab:high-rank}(c) we see how the TT
implementation failed to compute an approximation of $\us$ because the
cross-interpolation was unable to approximate $\us$ on a coarse grid
with a high rank.
The smaller speed-up observed with higher ranks compared to the lower
rank case can be explained by the fact that operations with higher
rank tensors require more floating point operations, partially
offsetting the advantage of using the tensor-train format.
%

\subsection{Discussion}
\label{subsec:numerical:discussion}

The performance of the previous sections demonstrates that the TT
approach offers a robust and highly efficient alternative to
traditional full-grid solvers while preserving the desired numerical
properties.
The usage of this type of representation leads to massive time and
memory savings while guaranteeing almost identical accuracy with
respect to the full-grid approach.
However, the implementation based on the tensor-train format implies
additional costs since it requires more sophisticated data structures
and a redesign of the numerical algorithms for handling operations in
the tensor format.
While this additional complexity of the PDE solvers pays off in terms
of efficiency, it makes the method's implementation more challenging
than the straightforward full-grid approach.
Based on the presented data and context, we can indeed identify some
critical points that can help understand when the TT approach can be
expected to perform better compared to the full-grid approach.
This issue will be explored in our future works.
However, it is worth noting that such critical issues are generally
outweighed by the significant advantages in computational efficiency
and storage requirements, particularly for larger problem sizes.

\begin{itemize}
\item \emph{Initialization Overhead.} On the coarser grids, the TT
  method shows similar or occasionally slightly higher computational
  times compared to the full-grid approach.  This suggests that some
  computational overhead occurs due to an initial TT decomposition and
  setup.
\item \emph{Memory Access Patterns.} The TT format, while more
  memory-efficient, may have less cache-friendly memory access
  patterns due to its smaller decomposed structure.  This could
  explain why the computational advantage is not as pronounced for
  smaller problem sizes where cache effects may be less significant.
\item \emph{Additional Parameters.} The TT method introduces
  additional parameters such as the accuracy threshold for the
  rounding after normal arithmetic operations or inside the conjugate
  gradient iterations, and the ranks of the initial decomposition.
\end{itemize}

Finally, a particularly noteworthy result emerges from these
calculations.
The explicit Euler method is traditionally considered inefficient for
parabolic problems due to the CFL stability condition, which requires
the time step size $\Delta t$ to scale with $\hh^2$.
This constraint means that when the computational grid is refined by
halving the spatial step size $\hh$, we must reduce $\Delta\ts$ by a
factor equal to $4$, thus implying a fourfold increase in the number
of iterations needed to reach a given final time $T$.
The stability constraint is less restrictive for the implicit Euler
and semi-implicit Crank-Nicolson methods, and it is a common practice
to halve $\Delta t$ when $\hh$ is halved, just to keep the balance
between the temporal and spatial approximation errors.
However, a comparison of the absolute time costs (i.e., columns
$\ctime_{TT}$) across all three time-marching schemes reveals an
unexpected advantage of the TT format.
Despite using an efficient preconditioned conjugate gradient method
for the implicit and semi-implicit schemes, the computational cost per
explicit Euler scheme iteration is much smaller in the TT format and
remains competitive with implicit methods.
This fact is in contrast with the usual (and accepted) notion
regarding the relative efficiency of explicit versus implicit methods
and deserve more investigation, which will be carried out in our
future works.


\section{Final Remarks}
\label{sec:conclusions}

In this paper, we have developed and analyzed a tensor-train low-rank
finite difference method for the efficient numerical solution of
three-dimensional parabolic problems.
Our approach combines the dimensional reduction capabilities of the
tensor-train format with classical finite difference discretization
techniques, resulting in a powerful computational framework that
significantly reduces both memory requirements and computational
complexity.

The key innovation of our method lies in the tensorization process
applied to traditional grid difference formulas, enabling the
representation of solutions in a low-rank tensor-train format.
This reformulation dramatically reduces the computational demands
while maintaining accuracy, thereby extending the practical
feasibility of solving three-dimensional problems that were previously
intractable due to the curse of dimensionality.
The integration of implicit and semi-implicit (Crank-Nicolson) time-marching schemes further enhances the method's robustness and
stability, providing a flexible framework which can be adapted to
diverse problem settings.

Our comprehensive numerical experiments demonstrate the method's superior performance compared to a conventional finite difference
method's design.
Specifically, studying how the algorithm performs in the case of a
rank-1 exact solution with a rank-1 initial approximation, we observe
that the algorithm is rank-stable in the sense that there is no
uncontrolled growth of the TT ranks.
Additionally, we observe a substantial reduction in memory
requirements while maintaining solution accuracy and a significant
improvement in computational efficiency, particularly for large-scale
problems.
The robust performance across various test cases and problem
parameters, also in the case of manufactured high-rank solutions, is promising for effectively handling challenging problem scenarios.

These results establish the practical viability of our approach for
applications in computational physics, engineering, and scientific computing.
The method's demonstrated effectiveness in handling three-dimensional parabolic problems suggests its potential applicability to a broader range of scientific and engineering challenges, particularly in areas where computational efficiency is crucial.
Future research directions may include the extension to a more general
classes of partial differential equations with application to specific
industrial and scientific problems.


\section*{Acknowledgments}
The LDRD-ER Project \#20210485ER and the LDRD-DR Project \#20250032DR
funded the work of Dr. G.~Manzini.
The LDRD-ER Project \#20210485ER also partially supported a three-week
visit of Dr. T.~Sorgente at LANL in August 2023.
Los Alamos National Laboratory is operated by Triad National Security,
LLC, for the National Nuclear Security Administration of the
U.S. Department of Energy (Contract No.\ 89233218CNA000001).
This work has been approved for public release and assigned the LA-UR number "LA-UR-25-22205".
T. Sorgente is a member of the RAISE Innovation Ecosystem, funded by the
European Union - NextGenerationEU and by the Ministry of University
and Research (MUR), National Recovery and Resilience Plan (NRRP),
Mission 4, Component 2, Investment 1.5, project ``RAISE - Robotics and
AI for Socio-economic Empowerment'' (ECS00000035).
G. Manzini and T. Sorgente are affiliated to the Italian Gruppo
Nazionale Calcolo Scientifico - Istituto Nazionale di Alta Matematica
(GNCS-INdAM).


\bibliographystyle{plain}


\clearpage


\newcommand{\Xsr}{\scalebox{0.725}{\Xs}}
\newcommand{\Ysr}{\scalebox{0.725}{\Ys}}
\newcommand{\Zsr}{\scalebox{0.725}{\Zs}}

\newcommand{\usX}{\us_{\Xsr}}
\newcommand{\usY}{\us_{\Ysr}}
\newcommand{\usZ}{\us_{\Zsr}}

\newcommand{\usXX}{\us_{\Xsr\Xsr}}
\newcommand{\usYY}{\us_{\Ysr\Ysr}}
\newcommand{\usZZ}{\us_{\Zsr\Zsr}}

\newcommand{\usXY}{\us_{\Xsr\Ysr}}
\newcommand{\usYZ}{\us_{\Ysr\Zsr}}
\newcommand{\usZX}{\us_{\Zsr\Xsr}}

\newcommand{\Dx}{\Delta\xs}
\newcommand{\Dy}{\Delta\ys}
\newcommand{\Dz}{\Delta\zs}

\appendix
\section{Taylor expansions for first derivatives}
\label{sec:appx:Taylor_expansion_1st_derivatives}

We recall that 
\begin{align*}
  \xsC(\ic) = \frac{ \xsV(\iV) + \xsV(\iV+1) }{2}\quad\textrm{and}\quad \hx = \xsV(\iV+1) - \xsV(\iV),\\[0.5em]
  \ysC(\jc) = \frac{ \ysV(\jV) + \ysV(\jV+1) }{2}\quad\textrm{and}\quad \hy = \ysV(\jV+1) - \ysV(\jV),\\[0.5em]
  \zsC(\kc) = \frac{ \zsV(\kV) + \zsV(\kV+1) }{2}\quad\textrm{and}\quad \hz = \zsV(\kV+1) - \zsV(\kV). 
\end{align*}
so that, for three given integers $m$, $n$, $p$ we can write 
\begin{align*}
  \begin{array}{rcll}
    \xsV(\iV)+\ms\hx &\!\!=\!\!& \xsC(\ic)+\Dx(m) &\quad\textrm{with}~\Dx(m) = (m-1/2)\hx,\\[0.5em]
    \ysV(\jV)+\ns\hy &\!\!=\!\!& \ysC(\jc)+\Dy(n) &\quad\textrm{with}~\Dy(n) = (n-1/2)\hy,\\[0.5em]
    \zsV(\kV)+\ps\hz &\!\!=\!\!& \zsC(\kc)+\Dz(p) &\quad\textrm{with}~\Dz(p) = (p-1/2)\hz.
  \end{array}
\end{align*}
From the definition of $\Dx(m)$, $\Dy(n)$, and $\Dz(p)$ we note that
$(\Dx(0))^2 = (\Dx(1))^2$, $\Dx(1)-\Dx(0) = \hx$, and
\begin{align*}
  \sum_{n=0}^{1}\Dy(n) = \sum_{p=0}^{1}\Dz(p) = 0.
\end{align*}
We are now ready to perform the Taylor expansion of the vertex grid
function values $\usV(\iV+\ms,\jV+\ns,\kV+\ps)$ around the center of
cell $\C(\ic,\jc,kc)$.
In such expansion, we denote the first and second derivatives of
$\us(x,y,z)$ with respect to $x$, $y$, $z$ at the cell center as
$\usX=\partial\us(\xsc,\ysc,\zsc)\slash{\partial\xs}$,
$\usXX=\partial^2\us(\xsc,\ysc,\zsc)\slash{\partial\xs^2}$, etc.
Since we assume that $\usV$ is obtained by sampling $\us(\xs,\ys,\zs)$
at the mesh vertices, we find that
\begin{align}
  &\usV(\iV+\ms,\jV+\ns,\kV+\ps)
  = \us( \xsV(\iV)+\ms\hx, \ysV(\jV)+\ns\hy, \zsV(\kV)+\ps\hz )       \nonumber\\[0.5em]
  &\qquad= \us( \xsC(\ic)+\Dx(m), \ysC(\jc)+\Dy(n), \zsC(\kc)+\Dz(p) )\nonumber\\[0.5em]
  &\qquad= \usC + \usX\Dx(m) + \usY\Dy(n) + \usZ\Dz(p)                \nonumber\\[0.5em]
  &\qquad\quad+ \frac12\Big( \usXX(\Dx(m))^2 + \usYY(\Dy(n))^2 + \usZZ(\Dz(p))^2 \Big)\nonumber\\[0.5em]
  &\qquad\quad+ \frac12\Big( \usXY\Dx(m)\Dy(n) + \usYZ\Dy(n)\Dz(p) + \usZX\Dz(p)\Dx(m) \Big)
  + \calO(\hh^3).
  \label{eq:Taylor:expansion}
\end{align}
Then, we sum over $\ns\in\{0,1\}$, and $\ms\in\{0,1\}$ and we find that
\begin{align*}
  \phi(m) := \sum_{n=0}^{1}\sum_{p=0}^{1} \usV(\iV+\ms,\jV+\ns,\kV+\ps)
  = 4\usC + 4\usX\Dx(m) + \big[\textrm{second~order~terms}\big] + \calO(\hh^3),
\end{align*}
where $\phi(m)$ is a suitably defined auxiliary function that depends
on $\ms\in\{0,1\}$, and, for simplicity, we collect all second-order
terms in the squared bracket term $[\,\ldots]\,$, which also depends on $\ms$.
It is easy to see that this latter term is the same for $m=0$ and
$m=1$.
Finally, we take the difference between the values of $\phi(m)$ for
$m=1$ and $m=0$, and using the definition of $\usVx$ we find that
\begin{align*}
  4\hx\usVx = \phi(1) - \phi(0) = 4\hx\usX + \calO(\hh^3),
\end{align*}
which implies the first assertion in
Eq.~\eqref{eq:2nd-order:accuracy}.
We can prove the other two relations in~\eqref{eq:2nd-order:accuracy} by repeating the same
argument.

\section{Taylor expansions for second derivatives}
\label{sec:appx:Taylor_expansion_2nd_derivatives}

We substitute the Taylor expansion~\eqref{eq:Taylor:expansion}
in~\eqref{eq:second:X}, which we rewrite here for convenience:
\begin{align*}
  \usVxx(\ic,\jc,\kc)
  = \frac{1}{8\hx^2}\sum_{m=0}^{3}\eta_{m}\,\sum_{n=0}^{1}\sum_{p=0}^{1} \xi_{n}\xi_{p} \usV(\iV-1+m,\jV+n,\kV+p),
\end{align*}
where we recall that the coefficients $\eta_{m}$ and $\xi_{n}$ (or
$\xi_{p}$) take the values $\xi_{0}=-1$, $\xi_{1}=1$,
$\eta_{0}=\eta_{3}=1$, $\eta_{1}=\eta_{2}=-1$.
Then, we note that
\begin{itemize}[nosep]
\item[-] the terms with $\us(\xsc,\ysc,\zsc)$ cancel out due to
  $\sum_{n=0}^{1}\sum_{p=0}^{1} \xi_{n}\xi_{p} = 0$;
\item[-] the first derivatives in $\ys$ and $\zs$ cancel out due to the
  symmetry of $\xi_{n}\xi_{p}$;
\item[-] the mixed derivatives cancel out for the same reason;
\item[-] the second derivatives in $\ys$ and $\zs$ only contribute to
  the generic error term $ \calO(\hh^2)$.
\end{itemize}
Therefore, we find that
\begin{align*}
  \usVxx = \usXX + \calO(\hh^2),
\end{align*}
which prove the first assertion
in~\eqref{eq:2nd-order:accuracy:2nd-derivatives}.
We can prove the other two relations
in~\eqref{eq:2nd-order:accuracy:2nd-derivatives} by repeating the same
argument.

\section{Grid transformation}
\label{sec:appx:grid_transformation}

We briefly sketch the derivation of Eq.~\eqref{eq:remapped:Laplacian}
to introduce the explicit transformation formulas for the first and
second derivatives, which will be useful in the derivation of the
finite difference formulas in both FG and TT formats.
Let $\ust:\widetilde{\Omega}\to\REAL$ and $\us:\Omega\to\mathbb{R}$ be
two functions related by the coordinate transformation, such that
$\ust(\xst)\equiv\us(\xs)$ whenever $\xst=\xst(\xs)$.
By applying the chain rule, we obtain the fundamental relation for the
first derivative:
\begin{align*}
  \frac{\partial}{\partial\xst}\Big(\ust(\xst)\Big)
  = \frac{\partial\us(\xs)}{\partial\xs}\,\frac{\partial\xs(\xst)}{\partial\xst}.
\end{align*}
This relation can be generalized to operate on any expression,
yielding:
\begin{align*}
  \frac{\partial}{\partial\xst}\Big[\cdot\Big] = 
  \left(\frac{\partial\xs(\xst)}{\partial\xst}\right)
  \frac{\partial}{\partial\xs}\Big[\cdot\Big].
\end{align*}
For second derivatives, applying the operator twice results in:
\begin{align*}
  \frac{\partial^2}{\partial\xst^2}\Big[\cdot\Big] =
  \left(\frac{\partial\xs(\xst)}{\partial\xst}\right)\frac{\partial}{\partial\xs}\,
  \left(\frac{\partial\xs(\xst)}{\partial\xst}\right)\frac{\partial}{\partial\xs}\Big[\cdot\Big].
\end{align*}
A crucial property of these transformations follows from their inverse
relationship.
Since $\xs(\xst(\xs))=\xs$, differentiating both sides yields:
\begin{align*}
  1 = \frac{\partial\xs(\xst(\xs))}{\partial\xs}
  = \frac{\partial\xs(\xst)}{\partial\xst}\frac{\partial\xst(\xs)}{\partial\xs},
\end{align*}
which implies the fundamental inverse relation:
\begin{align*}
  \frac{\partial\xs(\xst)}{\partial\xst}=\left[\frac{\partial\xst(\xs)}{\partial\xs}\right]^{-1}.
\end{align*}
Substituting this result into the second derivative expression, we obtain our final form:
\begin{align*}
  \frac{\partial^2}{\partial\xst^2}\Big[\cdot\Big] =
  \left(\frac{\partial\xst(\xs)}{\partial\xs}\right)^{-1}\frac{\partial}{\partial\xs}\,
  \left(\frac{\partial\xst(\xs)}{\partial\xs}\right)^{-1}\frac{\partial}{\partial\xs}\Big[\cdot\Big].
\end{align*}
This formulation is particularly useful when the inverse
transformation $\xst(\xs)$ is known explicitly, as it allows us to
express second derivatives in the physical domain in terms of
derivatives in the reference domain.


\section{Time integration: implementation details}
\label{sec:appx:time_integration}
Algorithm~\ref{algo:explicit:TT} details a single time step of the
explicit Euler method from $\ts^{n}$ to $\ts^{n+1}$.
The input of the algorithm consists in $\usV^{\TT,n}$, the TT grid
function at the time step $\ts^{n}$, the right-hand side function
$\fs^{\TT,n}$, and the boundary function $\gs^{n}$, respectively
evaluated at the internal and boundary vertices.
The output consists in the updated TT vertex grid function
$\usV^{\TT,n+1}$ at the next time step.
The algorithm first enforces the boundary conditions.
For all vertices residing in the ghost boundary frames described in
Section~\ref{subsec:FVM-regular:BCS}, the auxiliary vertex grid function
$\usVb^{\TT,n}$ is assigned $\gs(\xv,\ts^{n})$, the value of the
boundary function evaluated at the vertex spatial coordinates $\xsV$
and current time $\ts^{n}$.
Then, the algorithm stores the internal values of $\usV^{\TT,n}$ into 
the auxiliary field $\usVb^{\TT,n}$.
This preserves the solution at time $t^n$ before the field
$\usV^{\TT,n}$ is overwritten with the updated field $\usV^{\TT,n+1}$.
The algorithm computes the Laplacian of $\usVb^{\TT,n}$ at cell
centers using the discretization formulas introduced in
Section\ref{sec:FVM-regular}.
This operation is denoted as $\text{Compute Laplacian}(\usVb)$.
Then, the algorithm interpolates the cell-centered Laplacian
$\DeltaCTT\usVTTb$ at the mesh vertices resulting in
$\DeltaVTT\usVTTb$.
This step is necessary to obtain the Laplacian values at the same
locations as the vertex grid function.
Then, we compute $\Delta\ts(\DeltaVTT\usVTTb + \fsVTT)^{n}$, the
residual at each vertex, where $\fsV^{\TT,n}$ is the TT representation
of the right-hand side function evaluated at the mesh vertices and
current time $\ts^n$, as detailed in Section~\ref{sec:time_integration}.
The vertex grid function is then updated by overwriting
$\usVb^{\TT,n}$ with the new field $\usVb^{\TT,n+1}$.
Finally, the algorithm performs the rounding operation to control the
rank growth and assign the updated field to $\usV^{\TT,n}$.

\begin{algorithm}[hbt!]
  \caption{Explicit Euler Method}
  \label{algo:explicit:TT}
  \textbf{Input:}
  \begin{itemize}[nosep]
  \item $\fs$:          right-hand side function on $\Omega$
  \item $\gs$:          boundary function on $\partial\Omega$
  \item $\usV^{\TT,n}$:  TT vertex grid function at time $t^n$
  \end{itemize}
  \textbf{Output:}
  \begin{itemize}[nosep]
  \item $\usV^{\TT,n+1}$: TT vertex grid function at time $t^{n+1}$
  \end{itemize}
  \begin{algorithmic}[1]
    \Procedure{Explicit Time Marching Scheme}{$\fs, \gs, \usV^n$}
    \State $\usVb^{\TT,n}(\iV,\jV,\kV) \gets \gs(\xs(\iV),\ys(\jV),\zs(\kV),\ts^n) \quad \forall \mbox{\textrm{~boundary~vertex~}} \V(\iV,\jV,\kV)$  \hfill \Comment{Assign BCs}
    \State $\usVb^{\TT,n}(\iV,\jV,\kV) \gets \usz{\V}{\TT,n}(\iV,\jV,\kV)          \quad \forall \mbox{\textrm{~internal~vertex~}}\V(\iV,\jV,\kV)$   \hfill \Comment{Store Previous Values}
    \State $\DeltaCTT\usVTTb \gets \text{Compute Laplacian}(\usVTTb)$                                                                             \hfill \Comment{Compute Cell-Centered Laplacian}
    \State $\DeltaVTT\usVTTb \gets \text{Vertex Interpolation}(\DeltaCTT\usVTTb)$                                                                 \hfill \Comment{Interpolate Laplacian to Vertices}
    \State $\fsV^{\TT,n} \gets \gs$ \hfill \Comment{Compute forcing term}
    \State $\usVb^{\TT,n+1} \gets\usVb^{\TT,n} +\Delta\ts\left(\DeltaVTT\usVb^{\TT,n} +\fsV^{\TT,n} \right)\quad\forall\mbox{\textrm{~internal~vertex}}$  \hfill \Comment{Update internal vertices}
    \State $\usV^{\TT,n+1} \gets \RNDG(\usVb^{\TT,n+1})$                                                                                              \hfill \Comment{Apply rounding procedure}
    \EndProcedure
  \end{algorithmic}
\end{algorithm}

Algorithm~\ref{algo:implicit:TT} describes a single time step of the
implicit Euler method for evolving a solution $\usV^{\TT,n}$ at time
$\ts^n$ to $\usV^{\TT,n+1}$ at time $\ts^{n+1}$.
Like the explicit version, it takes the current solution
$\usV^{\TT,n}$, the right-hand side function $\fs$, and the boundary
condition function $\gs$ as inputs, and return the updated solution
$\usV^{\TT,n+1}$ at the following time instant $\ts^{n+1}$.
The auxiliary field in tensor-train format $\usVb^{\TT,n}$ is used to
collect the boundary values at the current time step and to store the
internal values of $\usV^{\TT,n}$ to preserve the solution at $\ts^n$.
The key difference from the explicit method is that the implicit Euler
method requires solving a linear system by applying the PCG to obtain
$\usV^{\TT,n+1}$.
This step is represented in Algorithm~\ref{algo:implicit:TT} by the
operation $\usV^{\TT,n+1} \gets \text{PCG}(\usVb^{\TT,n})$.

\begin{algorithm}[hbt!]
  \caption{Implicit Euler Method}
  \label{algo:implicit:TT}
  \textbf{Input:}
  \begin{itemize}[nosep]
  \item $\fs$:         right-hand side function on $\Omega$
  \item $\gs$:         boundary function on $\partial\Omega$
  \item $\usV^{\TT,n}$:  TT vertex grid function at time $t^n$
  \end{itemize}
  \textbf{Output:}
  \begin{itemize}[nosep]
  \item $\usV^{n+1}$: TT vertex grid function at time $t^{n+1}$
  \end{itemize}
  \begin{algorithmic}[1]
    \Procedure{Implicit Time Marching Scheme}{$\fs, \gs, \usV^n$}
    \State $\usVb^{\TT,n}(\iV,\jV,\kV) \gets \gs(\xs(\iV),\ys(\jV),\zs(\kV),\ts^n) \quad \forall \mbox{\textrm{~boundary~vertex~}} \V(\iV,\jV,\kV)$  \hfill \Comment{Assign BCs}
    \State $\usVb^{\TT,n}(\iV,\jV,\kV) \gets \usz{\V}{\TT,n}(\iV,\jV,\kV)          \quad \forall \mbox{\textrm{~internal~vertex~}}\V(\iV,\jV,\kV)$   \hfill \Comment{Store Previous Values}
    \State $\fsV^{\TT,n+1} \gets \gs$ \hfill \Comment{Compute forcing term}
    \State $\bsV=\RNDG(\usVb^{\TT,n}+\Delta\ts\fsV^{\TT,n+1})$ \hfill \Comment{Compute RHS term}
    \State $\usV^{n+1} \gets \text{PCG}(\bsV,\usVb^{\TT,n})$ \hfill \Comment{Solve Linear System using Preconditioned Conjugate Gradient}
    \EndProcedure
  \end{algorithmic}
\end{algorithm}

Algorithm~\ref{algo:Crank-Nicolson:TT} presents the semi-implicit
Crank-Nicolson method for evolving the solution $\usV^{\TT,n}$ at time
$\ts^n$ to $\usV^{\TT,n+1}$ at time $\ts^{n+1}$.
This algorithm leverages both explicit and implicit Euler steps to
achieve second-order accuracy in time.
The inputs remain the same as for the explicit and implicit Euler
methods: the current solution $\usV^{\TT,n}$, the right-hand side
function $\fs$, and the boundary condition function $\gs$.
The algorithm initially sets the boundary conditions for the explicit
step at the current time step $\ts^n$.
Then, an explicit Euler step is taken with a half time step
($\Delta\ts/2$) to compute an intermediate solution $\usV^{\TT,n+1/2}$
from the solution at $\ts^n$.
Before the implicit step, the boundary conditions are updated to the
values at the final time $\ts^{n+1}$.
Finally, an implicit Euler step, also with a half time step
($\Delta\ts/2$), is used to compute the solution $\usV^{\TT,n+1}$ from
the intermediate solution $\usV^{\TT,n+1/2}$.
As in the implicit Euler algorithm, this step involves solving a
linear system by applying the PCG method.
By combining the explicit and implicit steps in this manner, the
Crank-Nicolson method achieves second-order accuracy in time while
retaining some of the computational advantages of the explicit
method.
The use of half time step $\Delta\ts\slash{2}$ is essential for the
correct implementation of the Crank-Nicolson scheme.

\begin{algorithm}[hbt!]
  \caption{Semi-implicit Crank-Nicolson Method}
  \label{algo:Crank-Nicolson:TT}
  \textbf{Input:}
  \begin{itemize}[nosep]
  \item $\fs$:         right-hand side function on $\Omega$
  \item $\gs$:         boundary function on $\partial\Omega$
  \item $\usV^{\TT,n}$:  vertex grid function at time $t^n$
  \end{itemize}
  \textbf{Output:}
  \begin{itemize}[nosep]
  \item $\usV^{n+1}$: vertex grid function at time $t^{n+1}$
  \end{itemize}
  \begin{algorithmic}[1]
    \Procedure{Semi-implicit Time Marching Scheme}{$\fs, \gs, \usV^n$}
    \State \text{Set boundary conditions at $\ts=\ts^n$}
    \State \text{Compute $\usV^{\TT,n+1/2}$ from $\usV^{\TT,n}$ using the explicit Euler scheme with time step $\Delta\ts/2$}
    \State \text{Set boundary conditions at $\ts=\ts^{n+1/2}$}
    \State \text{Compute $\usV^{\TT,n+1}$ from $\usV^{\TT,n+1/2}$ using the implicit Euler scheme with time step $\Delta\ts/2$}
    \EndProcedure
  \end{algorithmic}
\end{algorithm}

\end{document}